\documentclass[11pt]{article}%
\usepackage{yfonts}
\usepackage{graphicx}
\usepackage{subcaption}
\usepackage{makeidx}
\usepackage{amssymb}
\usepackage[all]{xy}
\usepackage{float}
\usepackage[latin1]{inputenc}
\usepackage{amsfonts}
\usepackage{amsmath}
\usepackage{amssymb}
\usepackage{url}
\usepackage{hyperref}
\usepackage{graphicx}%
\usepackage{mathabx}
\usepackage{wrapfig}

\usepackage{amsthm}
\usepackage{url}
\usepackage{hyperref}
\providecommand{\U}[1]{\protect\rule{.1in}{.1in}}

\usepackage{xcolor}

\definecolor{ablue}{rgb}{0.36, 0.54, 0.66}
\definecolor{fuchsia}{rgb}{0.57, 0.36, 0.51}
\definecolor{agreen}{rgb}{0.0, 0.5, 0.0}
\definecolor{burn}{rgb}{0.43, 0.21, 0.1}
\definecolor{bgreen}{rgb}{0.0, 1.0, 0.0}

\newtheorem{prop}{Proposition}[section]
\newtheorem{cor}[prop]{Corollary}
\newtheorem{defi}[prop]{Definition}

\newtheorem{lem}[prop]{Lemma}

\newtheorem{theo}[prop]{Theorem}

\def\tr{\mbox{\rm Tr}}

\newcommand{\EE}{\mathbb{E}}

\newcommand{\NN}{\mathbb{N}}

\newcommand{\RR}{\mathbb{R}}

\newcommand{\Ba}{ {\cal B }}

\newcommand{\Na}{ {\cal N }}

\newcommand{\Ua}{ {\cal U }}

\newcommand{\Ga}{ {\cal G }}

\newcommand{\Ia}{ {\cal I }}
\newcommand{\Xa}{ {\cal X }}
\newcommand{\Ma}{ {\cal M }}

\newcommand{\Wa}{ {\cal W }}

\newcommand{\point}{\mbox{\LARGE .}}

\newcommand{\cqfd}{\hfill\blbx \\}
\def\blbx{\hbox{\vrule height 5pt width 5pt depth 0pt}\medskip}

\def \RR{\mathbb{R}}
\def \SS{\mathbb{S}}
\def \EE{\mathbb{E}}

\begin{document}

  \title{Matrix product moments in normal variables}
  \author{A. N. Bishop, P. Del Moral}


\maketitle

\begin{abstract}
Let $\Xa=XX^{\prime}$ be a random matrix associated with a  centered
$r$-column centered Gaussian vector $X$ with a covariance matrix $P$. In this article we compute expectations
of matrix-products of the form $\prod_{1\leq i\leq n}(\Xa P^{v_i})$ for any $n\geq 1$ and any multi-index
parameters $v_i\in\NN$. We derive closed form formulae and a simple sequential algorithm to compute these matrices w.r.t. the parameter $n$.
The second part of the article is dedicated to a non commutative binomial formula for the central matrix-moments $\EE\left(\left[\Xa-P\right]^n\right)$. The matrix product moments discussed in this study are expressed in terms of polynomial formulae w.r.t. the powers of the covariance matrix, with coefficients depending on the trace of these matrices. We also derive a series of estimates w.r.t. the Loewner order on quadratic forms. For instance we shall prove the rather crude estimate $\EE\left(\left[\Xa-P\right]^n\right)\leq \EE\left(\Xa^n-P^n\right)$, for any $n\geq 1$.\\

\newline


\end{abstract}



\section{Introduction}

Let $\SS_r$ be the set of symmetric $(r\times r)$-square matrices equipped with the Frobenius inner product 
$(S_1,S_2)\in \SS_r^2\mapsto \langle S_1,S_2\rangle_F=\tr(S_1S_2)$. We let $I$ be the identity matrix. Let $\SS_r^+\subset \SS_r$
be the cone of positive definite matrices equipped with the Loewner partial order $S_1\leq S_2$
between quadratic forms induced by the matrices $S_1,S_2\in \SS_r^+$.

We consider an $r$-column centered Gaussian vector $X$ with covariance matrix $P\in \SS_r^+$, and
set $\Xa:=XX^{\prime}$.  
When $P=I=\EE(\Xa)$, for any $n\geq 1$ we have the more or less well known moment formulae
\begin{eqnarray}
\EE(\Xa^{n})&=&\EE\left(\,\Vert X\Vert^{2(n-1)}~\Xa~\right)
\\
&=&2^{n-1}~(n-1)!
\sum_{0\leq k< n}~~\frac{2^{-k}}{k!}~
\EE\left[\Vert X\Vert^{2k} \right]~I\quad\mbox{\rm with}\quad \EE\left[\Vert X\Vert^{2k} \right]=\prod_{0\leq l<k}\left(r+2l\right)\nonumber
\end{eqnarray}
A proof of the above moment formula is provided in (\ref{Xa-X}).
In addition, we have the binomial formula
\begin{equation}\label{case-PI-intro}
\EE\left[\left(\Xa-I\right)^{n}\right]= \displaystyle  \frac{1}{r}
 \sum_{0\leq  k\leq n}
\left(
\begin{array}{c}
n\\
k
\end{array}
\right)~(-1)^{n-k}~\EE\left[\Vert X\Vert^{2k} \right]~I
\end{equation}
In the above display $\Vert \point\Vert$ stands for the Euclidian norm on $\RR^r$. When $r=1$ the above formula resumes to
$$
\EE\left[\left(X^2-I\right)^{n}\right]=
 \sum_{0\leq  k\leq n}
\left(
\begin{array}{c}
n\\
k
\end{array}
\right)~(-1)^{n-k}~\EE\left[X^{2k}\right]
$$
After some elementary computations we find that
\begin{equation}\label{case-PI-intro-estimate-r1}
0\leq \EE\left[(X^2-1)^{2m}\right]\leq ~\EE\left[X^{4m}\right]
-  \frac{3}{4} ~\frac{m-1/2}{m-3/4}~\left(1-\frac{2}{3m}\right)~
\sum_{1\leq k\leq m}~{2m \choose 2k-1}~\EE\left[X^{4k-2}\right]
\end{equation}
In addition, for any $r\geq 2$ we have the estimate
\begin{equation}\label{case-PI-intro-estimate-r2}
\begin{array}{l}
   \displaystyle\EE\left[\left(\Xa-I\right)^{2m}\right]
\leq \displaystyle  \EE\left[\Xa^{2m}\right]-
\frac{m}{2m+(r-2)}~\left(\frac{3}{2}+\frac{r-2}{2m-1}\right)
\sum_{1\leq k\leq m}~
\left(
\begin{array}{c}
2m\\
2k-1
\end{array}
\right)~\EE\left[\Xa^{2k-1}\right]
\end{array}
\end{equation}
Similar estimates can be derived for odd central moments. For completeness the detailed proof of the above formula
for even or odd powers are housed in the appendix on page~\pageref{proof-PI-intro-estimate-r12}. 

When $P=I$ we have the well known decomposition 
$\Xa= R~\Ua$ with $R=\langle X,X\rangle$, $\Ua=UU^{\prime}$, and $U=R^{-1/2}X$. It is well known that
$R$ and $U$ are independent. In addition we have $\tr(\Ua)=1$ with $\Ua\geq 0$ and $\Ua^2=\Ua$. The rank one random projections can be interpreted as density matrices. In quantum mechanics, these matrices are the quantum  version of probability densities. The computation of expectations trace of
powers of density matrices and related quantities is an active research area, see for instance~\cite{carlen} and references therein.

This article is mainly concerned with the extension of the formulae stated above to any covariance matrices. 

The random matrix $\Xa$ discussed above has a Wishart distribution with a single degree of freedom. 
Wishart distributions arise in a variety of fields and more particularly in statistics,
quantum physics, engineering, econometrics,  as well as in optimization theory, see for instance~\cite{bai,bao,caianiello,carlen,magnus,so,
tse-zeitouni}
and references therein. The moments of products of quadratic forms associated to these
random positive definite matrices is rather well developed. We refer the reader to the seminal works
by J. Magnus~\cite{magnus}, G. Letac and H. Massam~\cite{letac1} (see for instance section 3), and the more recent article by
S. Matsumoto~\cite{matsumoto}. 

Very few articles are concerned with the computation of
moments of product of random matrices. One of the main reasons is that the 
computation of these moments requires the use of sophisticated combinatorial tools, such as
the Isserlis' theorem  (also called the
 Wick formula).  This theorem expresses the product of the entries of Wishart matrices
 in terms of all possible matching between the indices of the product entries.

 The computation of 
matrix-products moments using these tools leads to sophisticated  and computationally 
too difficult to handle combinatorial structures and expressions. For instance, to the best of our knowledge
the Isserlis' theorem has never been used to compute the central matrix-moments of Wishart matrices.
Some attempts have been made in~\cite{zhu}, but the central matrix moments formulae stated in lemma 4 in~\cite{zhu} are much coarser than expected and remains without a proof. 

Product of the entries of matrices can be also related in a natural way to tensor products and multilinear algebras. In this perspective the combinatorial
analysis of quadratic form moments based on Isserlis' theorem can be connected pictorially to wiring diagrams, see for instance section 5.4 in~\cite{kueng}. Another commonly used technique to handle the combinatorial 
complexity of these expectations is the Weingarten Calculus~\cite{weingarten}.  This elegant combinatorial 
calculus allows to compute the moments of the Haar measure  on matrix algebras such as the law of the random matrices $\Ua$ discussed above. The central object is the Weingarten function on the symmetric group that encapsulates the combinatorial structure
of Wick formula. In this context all moments formulae are expressed in terms of sophisticated formulae based on Mobius inversion and convolutions
w.r.t. these functions. We refer the reader to the pioneering article by D. Weingarten~\cite{weingarten}, the article by C.B. Collins~\cite{collins}, and the recent article by 
B. Collins, S. Matsumoto~\cite{collins}. One of the main drawbacks of this integration calculus of matrix product elements is the combinatorial complexity of moments formulae expressed in terms of summations over the set of permutations.

 As mentioned above the first aim of this article is to compute expectations
of matrix-product of the form \begin{equation}\label{def-Ma}
\Ma(v):=(\Xa P^{v_n})\ldots (\Xa P^{v_n})
\end{equation} for any $n\geq 1$,  any multi-index
parameters $v_i\in\NN$, and any covariance matrix.
We derive a simple sequential algorithm and closed form formulae to compute these matrices w.r.t. the parameter $n$.
The second part of the article is dedicated to a non commutative binomial formula for the central matrix-moments $\EE\left(\left[\Xa-P\right]^n\right)$. For instance we shall prove the rather crude estimate $\EE\left(\left[\Xa-P\right]^n\right)\leq \EE\left(\Xa^n-P^n\right)$, for any $n\geq 1$.

Despite our efforts we have not succeeded in extending these results to more general
Wishart matrices. This question remains an important and open research problem that we hope to address in the near future. More detailed comments with a series of comparisons of our results with 
existing literature are provided in section~\ref{sec-comparison}.

\subsection{Statement of the main results}

We let $\Na:=\cup_{m\geq 0}\NN^m$ with the convention $\NN^0=\emptyset$. 
For any $v=(v_1,\ldots,v_m)\in \NN^m$ we set $\vert v\vert:=\sum_{1\leq i\leq m}v_i$.
For any $m\geq 1$ we also let $0_m\in\NN^m$ the $m$-row vector  of null entries. We consider the
positive map $M~:~\Na\rightarrow \SS_r^+ $ defined for any $m\geq 0$ 
by the matrix product moments
$$
M~:~v=(v_1,\ldots,v_m)\in\NN^m\mapsto M(v)=\EE\left(\left(\Xa P^{v_1}\right)
\cdots \left(\Xa P^{v_m}\right)\right)\in \SS_r^+
$$
with the convention $M(\emptyset)=I$ when $m=0$. The semi-definite positiveness property
$M(v)\in \SS_r^+$ is a consequence of matrix-polynomial formula (\ref{S-formula}) stated below.

We also consider the weighted matrix product moments
$$
\forall m\geq 1\quad \forall n\in\NN\qquad  W(m,n):=\sum_{}~(1+v_m)~M(v)\in \SS_r^+
$$
where the summation runs over all $v=(v_1,\ldots,v_m)\in \NN^m$ s.t. $\vert v\vert=n$
For instance when $n=1=m$  and $v\in \NN$ we have
$
M(v)=P^{1+v}\quad\mbox{\rm and}\quad W(1,1)=2~P^2
$. The weighted moments $W$ can be thought as the non commutative version of the binomial
coefficients, in the sense that 
\begin{equation}\label{case-PI}
P=I\Longrightarrow  W(m,n):= \left(
  \begin{array}{c}
  n+m\\
  n
  \end{array}
  \right)~
\EE(\Xa^{m})
\end{equation}

The last assertion comes from the fact that
\begin{equation}\label{C-n-m-n}
\displaystyle \sum_{v_1+\ldots+v_m=n}(v_{m}+1)= \left(
  \begin{array}{c}
  n+m\\
  n
  \end{array}
  \right)
\end{equation}
The proof of (\ref{C-n-m-n}) follows standard computations on binomial coefficients.
For the convenience of the reader some properties of the binomial coefficients and the detailed
proof of  (\ref{C-n-m-n})  is provided in the appendix on page~\pageref{proof-C-n-m-n}.

 The computation of the matrix product moments for general multi-indices
 requires lengthy and sophisticated calculations.

In this notation, our first main result takes basically the following form.
\begin{theo}\label{theo-1}
For any $(v,p)\in (\NN^m\times\NN)$ such that $\vert v\vert=n\geq 0$ and $m\geq 1$ we have the recursion
\begin{equation}\label{recursion-XX-moments-intro}
M(v,p)=P^{p+1}~\left[\tr(M(v))~I+2\sum_{1\leq i\leq m}~
M\left(v_1,\ldots,v_{i-1},v_{i+1},\ldots,v_m,v_i\right)\right]
\end{equation}
and the matrix-polynomial formula
\begin{equation}\label{S-formula}
M(v,p)=\sum_{0\leq q\leq m+n}~\rho^M_{m,n}(v,q)~P^{1+p+q}\quad\mbox{and}\quad
 W(m+1,n)=\sum_{0\leq q\leq m+n}~\rho^W_{m,n}(q)~P^{1+q}
\end{equation}
for some non negative parameters $\rho^M_{m,n}(v,q)$, resp.  $\rho^W_{m,n}(q)$, which can be explicitly expressed in terms of 
the parameters $(m,n,v,p)$, resp. $(m,n,v)$, and the
traces of  the covariance matrix powers. In addition, For any $0\leq q\leq m+n$ with $n\geq 0$ and $m\geq 1$ we have the monotone properties
 \begin{equation}
\rho^W_{m,n}(q) \geq 2\left(1-\frac{1}{n+2}\right)\rho^W_{m-1,n+1}(q) ~\mbox{and}~ 
W(m+1,n)\geq 2\left(1-\frac{1}{n+2}\right)W(m,n+1)
 \label{rhoW-inq}
\end{equation}
\end{theo}

The detailed descriptions of  the parameters $\rho^M_{m,n}(v,q)$ and $\rho^W_{m,n}(q)$ are  provided in section~\ref{some-notation-sec}, on page~\pageref{defi-rho}. We already mention that these formula are more or less direct consequences of the recursion (\ref{recursion-XX-moments-intro}).

The proof of the recursion  (\ref{recursion-XX-moments-intro}) is provided in the appendix on page~\pageref{recursion-XX-moments-intro-proof}. The proof of (\ref{rhoW-inq}) is  provided on page~\pageref{proof-rhoW-inq}. The proof of the l.h.s. formula in  (\ref{S-formula}) is provided on page~\pageref{S-formula-proof},
and the r.h.s. formula in (\ref{S-formula})  is proved in section~\ref{proof-sec-weighted}

The recursion (\ref{recursion-XX-moments-intro}) provides a simple but rather lengthy way to compute sequentially
the matrix polynomials $ M(v)$ 
without calculating sums over the set of permutations at every step. For instance, for any $(v_1,v_2,v_3)\in \NN^3$ we have
\begin{eqnarray*}
 M(v_1,v_2)&=& P^{1+v_{2}}~\left[
  ~\tr( M(v_1))~I~+~2~ M(v_1)\right] =~P^{1+v_{2}}~\left[
  ~\tr(P^{1+v_1})~I~+~2~P^{1+v_1}\right]
\end{eqnarray*}
and
\begin{eqnarray*}
 M(v_1,v_2,v_3)&=&~P^{1+v_{3}}~\left[ \tr(M(v_1,v_2))+2\left(M(v_2,v_1)+M(v_1,v_2)\right)\right]\\
 &=&P^{1+v_{3}}~\left[~\tr(P^{1+v_1})~\tr(P^{1+v_2})~I+2~\tr(
P^{2+v_1+v_2})\right.\\
&&\hskip1cm\left.+2 \left(\tr(P^{1+v_1})~P^{1+v_2}+\tr(P^{1+v_2})~
P^{1+v_1}\right)+2^3~P^{2+v_1+v_2}\right]
\end{eqnarray*}

Our second main result takes basically the following form.
\begin{theo}\label{theo-2}
For any $n\geq 0$ have the non commutative binomial formulae
\begin{equation}\label{Binomial-W}
\EE\left[\left(\Xa-P\right)^{n}\right]=(-1)^{n}~P^{n}+\sum_{0\leq k< n}~(-1)^k~W(n-k,k)
\end{equation}
In addition, we have the estimates
\begin{eqnarray}
 \EE\left[\left(\Xa-P\right)^{2n}\right]&\leq& \EE\left[\Xa^{2n}\right]-(2n-1)~P^{2n}-\frac{1}{4}~
 \sum_{1\leq k< n}~W(2(n-k)+1,2k-1)\nonumber\\
\EE\left[\left(\Xa-P\right)^{2n+1}\right]&\leq& \EE\left[\Xa^{2n+1}\right]-P^{2n+1} -\frac{1}{4}~\sum_{1\leq k\leq  n}~W(2(n+1-k),2k-1)\nonumber\\
&&\label{estimate-XX-n}
\end{eqnarray}
\end{theo}

The proof of the central moment formula (\ref{Binomial-W}) is provided in section~\ref{proof-sec-ae-Binomial-W}. See also proposition~\ref{prop-ae-binomial} for the almost sure version of the binomial formulae (\ref{Binomial-W}).

The proof of the estimates (\ref{estimate-XX-n}) is provided at the end of section~\ref{sec-m-e}, on page~\pageref{proof-estimate-XX-n}. 

When $P=I$ using (\ref{case-PI}) we recover (\ref{case-PI-intro-estimate-r2}) up to some proportional factor; that is we have
$$
(\ref{estimate-XX-n})\Longrightarrow
\EE\left[\left(\Xa-I\right)^{2m}\right]\leq \EE\left[\Xa^{2m}\right]-(2m-1)~I-\frac{1}{4}~
 \sum_{1\leq k< m}~ \left(
  \begin{array}{c}
  2m\\
  2k-1
  \end{array}
  \right)~
\EE(\Xa^{2k-1})
$$

Theorem~\ref{theo-2} can be used to derived several Laplace estimates. For instance
for any $0\leq 2t\leq 1/\tr(P)$ we have the rather crude estimates
$$
 \EE\left[e^{t\, (\Xa-P)}\right]\leq \EE\left[e^{t\,\Xa}\right]\quad \mbox{\rm and}\quad
  \partial_t~\EE\left[e^{t\, (\Xa-P)}\right]\leq \partial_t~\EE\left[e^{t\, \Xa}\right]
$$
After some elementary manipulations we also have 
$$
 \EE\left[e^{t\, (\Xa-P)}\right]+\left[\frac{t^2}{2}~P^2+\left(1+\frac{3t}{4}~P\right)\left[\sinh{(tP)}-tP\right]\right]
\leq \EE\left[e^{t\,\Xa}\right]$$

\subsection{Some notation}\label{some-notation-sec}

We let $\langle n\rangle^{\langle m\rangle}$ be the set of all $(n)_{m}:=n!/(n-m)!$ one to one mappings from $[m]:=\{1,\ldots,m\}$
into $ [n]$, with $m\leq n$. 
Let $\Ga_n=\langle n\rangle^{\langle n\rangle}$ be the symmetric group of permutations over the set $[n]$.
For any $\sigma\in \Ga_n$,  $\vert \sigma\vert$ stands for the number of cycles in $\sigma$, and
 $C(\sigma)$ the set of cycles of $\sigma$. To have some canonical decomposition of permutations into cycles, we 
   cyclically rearrange each cycle so that it begins with its largest element and 
cycles are ordered in increasing order of their largest elements. For any $v\in \NN^n$ and any $a\in \langle n\rangle^{\langle m\rangle}$ with $m\leq n$ we
set
\begin{equation}\label{def-tvx}
v_a:=\left(v_{a(1)},\ldots, v_{a(m)}\right)\qquad
v(a)=\sum_{1\leq i\leq m}~v_{a(i)}\quad\mbox{\rm and}\quad t_v(\Xa)=\prod_{1\leq i\leq n} \tr\left( \Xa P^{v_i}\right)
\end{equation}

In this notation, by theorem 3.1 in~\cite{letac1} for any $v\in\NN^n$ we have the trace formula
\begin{eqnarray}
t_v:=\tr(M(v))=\EE\left(t_v(\Xa)\right)=\sum_{\sigma\in \Ga_n}~2^{n-\vert \sigma\vert}~\prod_{c\in C(\sigma)}~\tr\left(P^{\vert c\vert+v(c)}\right)\label{Letac-formula}
\end{eqnarray}

For any $m,n,p\in \NN$ we also set
\begin{eqnarray*}
\Delta_{m,n}&:=& \{0,\ldots,m\}\times \{0,\ldots,n\}\qquad \Delta_{m,n}(p):=\{(k,l)\in \Delta_{m,n}~:~k+l=p\}\\
V_{m,n}&=&\{v=(v_1,\ldots,v_m)\in \NN^m~:~\vert v\vert=n\}
\end{eqnarray*}
 The polynomial descriptions (\ref{S-formula}) of the matrix product moments $(M,W)$
are defined in terms of the array 
of traces  $t_{v}(k,l)$  defined for any $(v,(k,l))\in (V_{m,n}\times\Delta_{m,n})$
 by
$$
t_{v}(k,l):=(m)_{k}^{-1}~
\sum_{a\in A_v(k,l)}~t_{v_a}\quad\mbox{\rm with}\quad
 A_{v}(k,l):=\{a\in \langle m\rangle^{\langle k\rangle}~:~v(a)=l\}
$$
Observe that for any $v\in V_{m,n}$ we have
$$
l<\wedge_{1\leq i\leq m} v_i\Longrightarrow \forall 0\leq k\leq m\quad t_{v}(k,l)=0
$$
and
$$
 \forall 0\leq l<n\quad A_{v}(m,l)=\emptyset\qquad   A_{v}(m,n)=\Ga_m\quad \mbox{\rm and}\quad
t_{v}(m,l)=1_{l=n}~t_v
$$
When $l=0$ we have $ V_{m,0}=\{0_m\}$  thus we have
$$
t_{0_m}(k,0)=\EE\left(t_{0_k}(\Xa)\right)=\EE(\langle X,X\rangle^k)=\sum_{\sigma\in \Ga_k}~2^{k-\vert \sigma\vert}~\prod_{c\in C(\sigma)}~\tr\left(P^{\vert c\vert}\right)
$$

We consider the coordinate projection mappings $\pi,\pi^+$ on $\Na$ defined for any
$n\geq 1$ and any $(v_1,\ldots,v_n)\in \NN^{n}$ by
 $$
\pi(v)=(v_1,\ldots,v_{n-1})\in \NN^{n-1}\quad \mbox{\rm and}\quad \pi^+(v)=v_n\in\NN
  $$
We let $T(m,n)\in\SS_r^+$ and $W(m,n)\in\SS_r^+$ be the collection of weighted 
matrix moments defined for any
 $(m,n)\in (\NN^{\star}\times \NN)$ by
$$
T(m,n):=\sum_{v\in V_{m,n}}~M(v) \quad \mbox{\rm and}\quad
  W(m,n):=\sum_{v\in V_{m,n}}~(1+\pi^+(v))~M(v)
$$
We also consider their traces 
$$
\tau(m,n):=\tr\left( T(m,n)\right)\quad\mbox{and}\quad
\varpi(m,n):=\tr\left( W(m,n)\right)
$$
For instance for $m=1$ we have
$$
T(1,n)=M(n)=P^{1+n}\quad\mbox{\rm and}\quad
  W(1,n):=~(n+1)~
M(n)=(n+1)~P^{1+n}
$$

For any $m,n,p\in \NN$ and $v\in V_{m,n}$ we   consider the collection of trace functions
\begin{eqnarray}
\rho^M_{m,n}(v,q)&:=&\sum_{(k,l)\in \Delta_{m,n}(q)}~2^{k}~(m)_k~~t_{v}(m-k,n-l)\nonumber\\
\rho^W_{m,n}(q)&:=&\sum_{(k,l)\in \Delta_{m,n}(q)}
~2^{k}~(m)_k~~~\left(
\begin{array}{c}
q+1\\
l
\end{array}
\right)~~
\tau(m-k,n-l)\label{defi-rho}
\end{eqnarray}

We use the convention 
$ \tau(0,n)=1_{n=0}
$ so that for any $q<m+n$
\begin{equation}\label{convention-rho}
\rho^W_{m,n}(q)=\sum_{(k,l)\in \Delta_{m-1,n}(q)}
~2^{k}~(m)_k~~~\left(
\begin{array}{c}
q+1\\
l
\end{array}
\right)~~
\tau(m-k,n-l)
\end{equation}
For $q=m+n$ we have the boundary value
$$
\rho^W_{m,n}(m+n)=2^{m}~m!~~~\left(
\begin{array}{c}
m+n+1\\
n
\end{array}
\right)~
$$
Observe that
\begin{eqnarray*}
\rho^M_{m,0}(0_m,m-q)&=&\sum_{(k,l)\in \Delta_{m,0}(m-q)}~2^{k}~(m)_k~~
~\EE(\langle X,X\rangle^{m-k})\\
&=&2^{m-q}~(m)_{m-q}~~
\EE(\langle X,X\rangle^{q})\Longrightarrow \frac{2^{-m}}{m!}~\rho^M_{m,0}(0_m,m-q)=\frac{2^{-q}}{q!}~\EE(\langle X,X\rangle^{q})
\end{eqnarray*}
and
\begin{eqnarray*}
M(0_{m+1})&=&\EE(\Xa^{m+1})\\
&=&\sum_{0\leq q\leq m}~\rho^M_{m,0}(0_m,q)~P^{1+q}=
\sum_{0\leq q\leq m}~2^{q}~(m)_{q}~~
\EE(\langle X,X\rangle^{m-q})~P^{1+q}\\
&=&2^m~m!~\sum_{0\leq q\leq m}~\frac{2^{-q}}{q!}~
\EE(\langle X,X\rangle^{q})~P^{m+1-q}
\end{eqnarray*}
This yields the formula
\begin{equation}\label{Xa-X}
\frac{2^{-n}}{n!}~\EE(\Xa^{n})=\frac{2^{-n}}{n!}~M(0_n)=\frac{1}{2n}
\sum_{0\leq k< n}~~\frac{2^{-k}}{k!}~
\EE\left(\langle X,X\rangle^k\right)~~P^{n-k}
\end{equation}

\section{Some preliminary results}
\subsection{Complete Bell polynomials}

Taking the trace in (\ref{Xa-X}) we obtain for any $n\geq 1$ the recursion
\begin{eqnarray}
\frac{2^{-n}}{n!}~\EE\left(\langle X,X\rangle^n\right)&=&\frac{1}{2n}
\sum_{0\leq k< n}~~\frac{2^{-k}}{k!}~
\EE\left(\langle X,X\rangle^k\right)~~\tr(P^{n-k})\label{recursion-E-X}\\
&=&\frac{1}{2n}~\tr(P^{n})+\frac{1}{2n}
\sum_{1\leq k< n}~~\frac{2^{-k}}{k!}~
\EE\left(\langle X,X\rangle^k\right)~\tr(P^{n-k})\nonumber\\
&=&\frac{1}{n!}~\sum_{\sigma\in\Ga_n}\left(\frac{1}{2}\right)^{\vert \sigma\vert}~\prod_{c\in C(\sigma)}
\tr\left(P^{\vert c\vert}\right)\nonumber
\end{eqnarray}
The last assertion is a direct consequence of the trace formula (\ref{Letac-formula}) obtained by 
Letac and Massam in~\cite{letac1} (see for instance theorem~\ref{letac-theo} in the present article). Formula (\ref{recursion-E-X}) can also be checked using the cycle index formula of the symmetric group. 
Recalling that $\Xa$ is distributed according to the  Wishart distribution
with a single degree of freedom and covariance matrix $P$, for any $2t<\tr(P)$ we also have
\begin{eqnarray}
\displaystyle\EE\left(\exp{\left(t\,\langle X,X\rangle\right)}\right)
&=&\exp{\left[\frac{1}{2}\sum_{n\geq 1}~\frac{(2t)^n}{n}~\tr(P^n)\right]}\label{t-X-Laplace}\\
\displaystyle\EE\left(\Xa~\exp{\left(t\Xa\right)}\right)
&=&\EE\left[\exp{\left(t\,\langle X,X\rangle\right)}\right]~[I-2t~P]^{-1}P
\label{alternative-recursion}
\end{eqnarray}
 Thus, the formula (\ref{Xa-X}) and the recursion (\ref{recursion-E-X})  can also be deduced from 
the Laplace formulae (\ref{t-X-Laplace}) and (\ref{alternative-recursion}). The proof of (\ref{t-X-Laplace}) and  (\ref{alternative-recursion})
follows standard calculations. To be completed the  detailed proofs are 
provided in the appendix, on page~\pageref{proof-estimate-Laplace-xx}.

Next we provide
a simple alternative proof and an explicit formula in terms of complete Bell polynomials.

The exponential Faa di Bruno formula applied to (\ref{t-X-Laplace}) reads
\begin{eqnarray*}
\EE\left(e^{t\langle X,X\rangle}\right)
&=&\sum_{n\geq 0}~\frac{t^n}{n!}~\Ba_n(b_1(P),\ldots,b_n(P))\quad\mbox{\rm with}\quad b_n(P):=\left[2^{n-1}(n-1)!~\tr(P^n)\right]
\end{eqnarray*}
with the complete Bell polynomials given by the formula
\begin{equation}\label{explicit-collins}
\begin{array}{l}
\Ba_{n}(b_1,\ldots,b_{n})
\displaystyle=n!~\sum_{1\leq p\leq n}~
\sum_{k_0:=0< k_1<k_2\ldots<k_{p-1}<k_{p}:=n}~\left\{
\prod_{l\leq l\leq p}\frac{1}{k_l}~\frac{b_{k_{l}-k_{l-1}}}{(k_l-(k_{l-1}+1))!}
\right\}\end{array}
\end{equation}
When $b_n=b_n(P)$ we find that
\begin{eqnarray*}
\frac{2^{-n}}{n!}~\EE\left(\langle X,X\rangle^n\right)&=&
\sum_{1\leq p\leq n}~~
\sum_{k_0:=0< k_1<\ldots<k_{p}:=n}~
\prod_{l\leq l\leq p}\frac{\tr(P^{k_{l}-k_{l-1}})}{2k_l}\end{eqnarray*}

The above formula can be proved directly from the recursion (\ref{recursion-E-X}).  For a more thorough discussion on Bell polynomial we refer to~\cite{collins}.

Combining this formula with (\ref{Xa-X}) we obtain the following corollary.
\begin{cor}\label{cor-1}
For any $n\geq 1$ we have
$$
\begin{array}{l}
\displaystyle~\frac{2^{-n}}{n!}~\EE(\Xa^{n})\\
\displaystyle=~\frac{1}{2n}~P^{n}+~\frac{1}{2n}
\sum_{1\leq  q<n}~\left[\sum_{1\leq p\leq n-q}~~
\sum_{k_0:=0< k_1<\ldots<k_{p}:=n-q}~
\prod_{l\leq l\leq p}\frac{\tr(P^{k_{l}-k_{l-1}})}{2k_l}\right]~~P^{q}
\end{array}
$$
\end{cor}
\begin{cor}
For any $n\geq 1$ we have
\begin{equation}\label{estimate-xx}
\frac{2^{-n}}{n!}~\EE(\Xa^{n}) \leq \displaystyle ~\frac{1}{2n}~P^{n}+~\frac{1}{2n}
\sum_{1\leq k<n}~\frac{(2k)!}{2^{2k}~k!^2}~\tr(P)^k~~P^{n-k}
\end{equation}

\end{cor}

The estimate (\ref{estimate-xx}) is a direct consequence of the following lemma of its own interest also proved  in the appendix, on page~\pageref{proof-lem-tech}. 
\begin{lem}\label{lem-tech}
We let $(f,g)$ be some non negative functions on $\NN$ such that 
$g(n)\in ]0,1/2]$ and for any $n\geq 1$ 
$$
f(n)\leq \frac{1}{n}\sum_{0\leq k<n}~g(n-k)~f(k)\quad\mbox{and}\quad \rho(g):=\sup_{n\geq 1}{\left[
g(n+1)/g(n)\right]}<\infty
$$
In this situation, for any $n\geq 1$ we have
\begin{eqnarray*}
f(n)&\leq& 
\rho(g)^{n}~\frac{2^{-2n}(2n)!}{n!^2}~
~f(0)\leq 
\rho(g)^{n}~\frac{1}{\sqrt{\pi n}}~\exp{\left(-\frac{1}{12n}\left[1-\frac{1}{9 n}\right]\right)}
~f(0)
\end{eqnarray*}
\end{lem}
Indeed, applying this lemma to the recursion (\ref{recursion-E-X}) with
$$
f(n)=\frac{2^{-n}}{n!}~\frac{\EE\left(\langle X,X\rangle^n\right)}{\tr(P)^n}\quad\mbox{\rm and}
\quad g(n)=\frac{1}{2}~\frac{\tr(P^{n})}{\tr(P)^{n}}\Longrightarrow \frac{g(n+1)}{g(n)}=
\frac{\tr(P^{n+1})}{\tr(P)\tr(P^{n})}\leq 1
$$
we find that
$$
\frac{2^{-n}}{n!}~\EE\left(\langle X,X\rangle^n\right)\leq \frac{2^{-2n}(2n)!}{n!^2}~\tr(P)^n
$$
The end of the proof is now a direct consequence of (\ref{Xa-X}). 
 
 \subsection{Extended binomial coefficients}
 
The extended binomial coefficient associated with some function $\beta~:n\in \NN\mapsto \beta_n\in \NN
$ are
given by the formula
 $$
 \left(
\begin{array}{c}
m+n-1\\
n
\end{array}
\right)_{\beta}:=
 \sum_{v\in V_{m,n}}~\beta_v~\quad\mbox{\rm with}\quad \beta_v:=\beta_{v_1}\ldots \beta_{v_{m}}
 $$
 with the convention $$ \left(\sum_{\emptyset},\prod_{\emptyset}\right)=(0,1)\Longrightarrow\left(
\begin{array}{c}
-1\\
0
\end{array}
\right)_{\beta}=1$$ when $m=n=0$. 
We recover the conventional definition of the binomial coefficients when $\beta_n=1$ is the unit function. 
For a more thorough discussion on the combinatorial properties of extended binomial coefficients and their relations with partial Bell polynomials
we refer the reader to the articles~\cite{cvi,fassi,eger}.

We set
$$
\Wa(n-k,k):=\sum_{0\leq k_1+k_2\leq k}
\left(
\begin{array}{c}
(n-1)-(k_1+k_2)-1\\
k-(k_1+k_2)
\end{array}
\right)_{t(\Xa)}
P^{k_1}\Xa P^{k_2}
$$
with the mapping $n\in \NN\mapsto t_n(\Xa)$ defined in (\ref{def-tvx}). We also have
\begin{equation}\label{case-PI-ae}
P=I\Longrightarrow  \Wa(m,n):= \left(
  \begin{array}{c}
  n+m\\
  n
  \end{array}
  \right)~
\Xa^{m}
\end{equation}
As in (\ref{case-PI}), we prove this claim using (\ref{C-n-m-n}) and the binomial formula stated in lemma~\ref{bino-appendix} in the appendix on page~\pageref{bino-appendix}.

\begin{prop}\label{prop-ae-binomial}
For any $n\geq 0$ have the almost sure non commutative binomial formulae
\begin{equation}\label{ae-Binomial-W}
\displaystyle\left(\Xa-P\right)^{n}=(-1)^{n}~P^{n}+\sum_{ 0\leq k<n}
~(-1)^{k}~\Wa(n-k,k)
\end{equation}
\end{prop}

The proof of (\ref{ae-Binomial-W}) is provided in section~\ref{proof-sec-ae-Binomial-W}.

Combining (\ref{case-PI-ae}) with (\ref{ae-Binomial-W}) when $P=I$ the binomial formula (\ref{ae-Binomial-W}) resumes to
\begin{equation}\label{case-PI-2}
\displaystyle\left(\Xa-I\right)^{n}=\sum_{ 0\leq k\leq n}
~(-1)^{k}~\left(
  \begin{array}{c}
  n\\
  k
  \end{array}
  \right)~\Xa^{n-k}
\end{equation}

\subsection{Comparisons with existing literature}\label{sec-comparison}

The literature in statistics abounds with formulae to compute expectations of product of quadratic forms
$\prod_{1\leq i\leq m}\langle X,Q_i X\rangle$ where $Q_i$ is a given collection of symmetric matrices.
See for instance~\cite{kendall,letac1,magnus,matsumoto}, and references therein. When all the matrices 
$Q_i=Q$ are equal by lemma 2.3 in~\cite{magnus} we have
$$
\EE\left(\langle X,Q X\rangle^m\right)=\EE\left(\langle \overline{X},\overline{X}\rangle^m\right)=\tr(\overline{M}(0_m))
$$
where $\overline{X}$ is an $r$-column centered Gaussian vector $X$ with a covariance matrix $\overline{P}:=P^{1/2}QP^{1/2}$, and $\overline{M}(0_m)=\EE\left(\overline{\Xa}^m\right)$ with $\overline{X}\,\overline{X}^{\prime}$.
In this situation, the recursions (\ref{recursion-XX-moments-intro}) and
 (\ref{recursion-E-X}), as well as corollary~\ref{cor-1}
can be used to compute sequentially or by summation over increasing sequences of integers these
$n$-moments. This result is rather well known. For instance by lemma 2.2 in Magnus~\cite{magnus}
 these moments are known to be polynomial in the variables
$\tr(P^k)$. The first $10$ moments are also given in Kendall and Stuart~\cite{kendall}.

As usual, the general case can deduced from the one discussed above using the polarization formulae
\begin{equation}\label{polarization}
n!~\EE\left(\prod_{1\leq i\leq n}\langle X,Q_i X\rangle \right)=\EE\left( \epsilon_W~\left\langle X,Q_WX\right\rangle^{\,n}\right)=\EE\left( \epsilon_W~\left\langle X_W,X_W\right\rangle^{\,n}\right)
\end{equation}
with
$$
 \epsilon_W:=\prod_{1\leq i\leq n}W_i
\qquad
Q_W:=\sum_{1\leq i\leq n}W_i~Q_i\in \SS_r\quad\mbox{\rm and}\quad X_W~\sim~\Na\left(0_r,P^{1/2}Q_WP^{1/2}\right)
$$
which is valid for any sequence of centered and independent random variables $(W_i)_{i\geq 1}$ with unit variance, and independent of $X$.
This polarization formula is rather well known. For completeness a detailed proof is provided in the appendix, on page~\pageref{proof-polarization}. The main drawback of this polarization technique is that it gives polynomial
formulae with signed coefficients. Thus it cannot be used to derive useful estimates for central moments.

As mentioned in the introduction, very few articles are concerned with computation of product matrix-moments.
Two notable exceptions can be underlined.
 
In~\cite{letac1} Letac and Massam design an elegant approach to
 compute  expectations of polynomials functions $P(\Xa)$ {\em w.r.t. the entries of}  Wishart random matrices  $\Xa$ with any degree of freedom,
 under the assumption that $P(\Xa)$ 
only depends on the eigenvalues of the matrix $\Xa$. The prototype of polynomials consider in this article are the determinant $\mbox{\rm det}(\Xa)$ and the trace of powers $\tr(\Xa^m)^n$. 
 As underlined by the authors these sophisticated mathematical objects doesn't provide a simple closed form solution but a relatively simple algorithm to compute expectations of 
power polynomial such as 
 $\Xa^n$ or $\Xa^{-1}$. 
 The article by 
 C. Kim and C. Kang~\cite{kim}
also proposes a tensor approach to compute the skewness of Wishart matrices. 

Using Isserlis' theorem Zhu proves in~\cite{zhu} that for any $n\geq 2$ we have
\begin{equation}\label{zhu-estimate}
\frac{2^{-n}}{n!}~\EE(\Xa^{n})\leq ~~
\frac{1}{8}~\tr(P)^{n-1}~P+~\frac{1}{4}~\tr(P)^{n-2}~P^2
\end{equation}
The author also states without proof that for any $n\geq 5$ we have
\begin{equation}\label{zhu-estimate-bis}
\EE((\Xa-P)^{n})\leq \EE(\Xa^{n})+P^n
\end{equation}
The estimates (\ref{estimate-XX-n}) are rigorously proved and they clearly improve (\ref{zhu-estimate-bis}).

Also observe that the estimate (\ref{estimate-xx}) yields the rather crude upper bound
\begin{eqnarray*}
 \frac{2^{-n}}{n!}~\EE(\Xa^{n})&\leq&
 \frac{1}{2n~\sqrt{\pi (n-1)}}~\tr(P)^{n-1}
~P\\
&&\hskip3cm
+\frac{1}{2n}~\left[P^{n-2}+
\sum_{3\leq k\leq n}~\frac{1}{\sqrt{\pi (k-2)}}~\tr(P)^{k-2}
~~P^{n-k}\right]~P^2\\
&\leq & \frac{1}{2n~\sqrt{\pi (n-1)}}~\tr(P)^{n-1}
~P
 +\frac{1}{n}\left[\frac{1}{2}+~\frac{1}{\sqrt{\pi }}
\left(\sqrt{n-1}-1\right)\right]~\tr(P)^{n-2}~P^2
\end{eqnarray*}
as soon as $n>1$.
In the last assertion we have used the fact that $P^k\leq \tr(P)^k~I$ and
$$
\sum_{1\leq k\leq n-2}~\frac{1}{2\sqrt{k}}\leq \sum_{1\leq k\leq n-2}~
\int_{k}^{k+1}\frac{1}{2\sqrt{t}}~dt=\sqrt{n-1}-1
$$
The above estimate improves (\ref{zhu-estimate}) as soon as $n\geq 3$.

 \section{Matrix product moments}
 \subsection{A matrix version of Letac-Massam trace formula}
 For any $n\geq 1$,
we equip the cartesian product $\SS_r^n$ with the product operations
 $$
 P\bullet Q=\left(P_{1}Q_{1},\ldots, P_{n}Q_{n}\right)\quad \mbox{\rm and the directed product}\quad
P^n= \prod_{1\leq i\leq n}P_i=P_1P_2\ldots P_n
 $$
for any $P=\left(P_{1},\ldots, P_{n}\right)$ and $Q=\left(Q_{1},\ldots, Q_{n}\right)$. When $P_i=P$
sometimes we write $P$ instead of $\left(P,\ldots, P\right)$. In this notation for any $P\in \SS_r$ we have
$$
 P\bullet Q=\left(PQ_{1},\ldots, PQ_{n}\right)\quad \mbox{\rm and}\quad P^n=\underbrace{P\ldots P}_{\mbox{\footnotesize $n$-times}}
$$
In this notation, for any $v=(v_1,\ldots,v_n)\in\NN^n$ we have
$$
P^v=\left(P^{v_1},\ldots,P^{v_n}\right)\in (\SS_r^+)^n\Longrightarrow  M(v)=\left(\EE\left[\Xa \bullet P^v\right]^n\right)\in \SS_r^+
$$

For any  $\sigma\in \Ga_n$ and $i\in [n]$ we let $c(i,\sigma)$ be the
 the cycle of $\sigma$ containing $i$ and set
 $C(i,\sigma):=C(\sigma)-\{c(i,\sigma)\}$.
 
 Let $
Q=\left(Q_{1},\ldots, Q_{n}\right)\in \SS_r
$ be a sequence of symmetric matrices. For a cycle $c$ of length $\vert c\vert=q$ we set
 $$
 Q_c:=\left(Q_{c_1},Q_{c_2},\ldots, Q_{c_q}\right) \qquad
 Q^c:=Q_{c_1}\ldots Q_{c_q}\quad\mbox{\rm and}\quad
 \tr_{B}\left(Q\right)=\prod_{c\in B}~\tr\left(Q^{c}\right)
 $$
 for any subset $B\subset C(\sigma)$ of cycles.
 
For any $v=(v_1,\ldots,v_n)\in\NN^n$ we also set $\overline{v}=(\overline{v}_1,\ldots,\overline{v}_n)\in\NN^n$
 with $\overline{v}_i=1+v_i$.
 
 \begin{defi}
We let $c^{\flat}(i,\sigma)$ be the {\em rooted} cycle 
deduced from $c(i,\sigma)$ by removing the state $i$
with adding the edge $ c_{p-1}(i,\sigma)\rightarrow c_{p+1}(i,\sigma)$. We also assume that 
the root of the cycle $c^{\flat}(i,\sigma)$ is
$c_{p+1}(i,\sigma)$; that is, we have
$$
c^{\flat}(i,\sigma):=(c_{p+1}(i,\sigma) \rightarrow  \ldots \rightarrow c_{q}(i,\sigma) \rightarrow c_{1}(i,\sigma) \rightarrow \ldots\rightarrow c_{p-1}(i,\sigma)\rightarrow c_{p+1}(i,\sigma))
$$
\end{defi}

\begin{prop}[Letac-Massam~\cite{letac1}]\label{letac-theo}
For any $n\geq 1$ and any
$
Q=\left(Q_{1},\ldots, Q_{n}\right)\in \SS_r
$
we have
$$
\EE\left(
\left[\Xa\bullet Q\right]^n\right)=2^n~\sum_{\sigma\in \Ga_n}
\left(\frac{1}{2}\right)^{\vert \sigma\vert}~\tr_{C(n,\sigma)}\left(P\bullet Q\right)~\left[\left(P\bullet Q\right)^{c^{\flat}(n,\sigma)}P\right]_{ sym}Q_n
$$
In particular for any $v\in\NN^n$ we have  the matrix polynomial formula
\begin{equation}\label{f1}
~ M(v)=2^n~\sum_{\sigma\in \Ga_n}
\left(\frac{1}{2}\right)^{\vert \sigma\vert}~\tr_{C(n,\sigma)}(P^{ \overline{v}})~P^{ \overline{v}(c(n,\sigma))}
\end{equation}

\end{prop}
The proof of this proposition is housed in the appendix on page~\pageref{proof-letac-theo}.

 \subsection{A polynomial expansion}
This section is mainly concerned with the proof of the l.h.s. formula  (\ref{S-formula}).

We associate with $v\in \NN^{m+1}$ and $i\in [m]$ the mappings
 $$
 \delta_{i}(v):=\left(\delta_{i}^{-}(v),v_{i}\right)\quad \mbox{\rm with}\quad
\delta_{i}^{-}(v) :=\left(v_1,\ldots,v_{i-1},v_{i+1},\ldots,v_{m}\right)
 $$
In this notation the recursion (\ref{recursion-XX-moments-intro}) takes the form
$$
\begin{array}{l}
\displaystyle  M(v)
\displaystyle=P^{\pi^+(v)+1}~\left[
\tr\left( M( \pi(v))\right)~I~+~2~\sum_{1\leq i\leq m}~
 M(\delta_{i}(v))\right]
\end{array}
$$
for any $m\geq 1$ and $v\in \NN^{m+1}$.
More generally, for any $1\leq k\leq m$  and any one to one mapping $a\in \langle m\rangle^{\langle k\rangle}$
we let 
$$
\delta_{a}^{-}~:~v\in \NN^{m+1}\mapsto \delta_{a}^{-}(v)\in \NN^{m-k}
$$ 
be the sequence defined as $v$ by deleting the $(m+1)$-th coordinate and ordered sequence of 
coordinates
$$
a(\sigma(1))<\ldots<a(\sigma(k))
$$
where $\sigma\in \Ga_k$ stands for an ordering permutation. More precisely, for any $v=(v_1,\ldots,v_{m+1})$ 
we have
$$
\delta_{a}^{-}(v):=\left(v_1,\ldots,v_{a(\sigma(1))-1},v_{a(\sigma(1))+1},\ldots,v_{a(\sigma(k))-1},v_{a(\sigma(k))+1},\ldots, v_{m}\right)\in \NN^{m-k}
$$
\begin{prop}
For any $m\geq 1$ and $v\in \NN^{m+1}$ we have the matrix polynomial formula
\begin{equation}\label{explicit-moments}
\begin{array}{l}
\displaystyle  M(v)-2^{m}~m!~~
P^{\vert\overline{v}\vert}\\
\\
\displaystyle=
\displaystyle ~P^{\overline{v}_{m+1}}~\left[
\tr\left( M(\pi(v))\right)~I~+~\sum_{1\leq k<m}~2^k~
\sum_{a\in \langle m\rangle^{\langle k\rangle}}~
P^{\overline{v}(a)}~\tr\left( M(\delta_{a}^{-}(v))\right)\right]
\end{array}
\end{equation}

\end{prop}
\proof

For any $a\in \langle m\rangle^{\langle k\rangle}$ we set
 $$
 \delta_{a}~:~v\in \NN^{m+1}\mapsto  \delta_{a}(v)=\left(\delta_{a}^{-}(v),v_{a(k)}\right)\in \NN^{m+1-k}
 $$
 When $k=1$ sometimes we write $\left(\delta_{a(1)}^{-}, \delta_{a(1)}\right)$ instead of $\left(
 \delta_{a}^{-}, \delta_{a}\right)$.
 
 Observe that for any $v\in \NN^{m+1}$ and any $1\leq i<m-k$
the $i$-th coordinate of
$$
\begin{array}{l}
\displaystyle\delta_{a}(v)\\
\\
:=\left(v_1,v_2,\ldots,v_{a(\sigma(1))-1},v_{a(\sigma(1))+1},\ldots,v_{a(\sigma(k))-1},v_{a(\sigma(k))+1},\ldots, v_{m-1},v_{m},v_{a(k)}\right)\in \NN^{m+1-k}
\end{array}
$$
is given by an unique index $b(i)\in [m]-a[k]$; that is we have that
$$
\delta_{a}(v)_i=v_{b(i)}
$$
We let  $a_i$ be the extension of $a\in \langle n\rangle^{\langle k\rangle}$ to the set $[k+1]$ defined by
 $$
 a_i~:~j\in ~[k+1]\mapsto~a_i(j)=1_{[k]}(j)~a(j)+1_{j=k+1}~b(i)
 $$
 Observe that $$
 a \in \langle m\rangle^{\langle k\rangle}\mapsto a_i\in  \langle m\rangle^{\langle k+1\rangle}
 $$
and
\begin{equation}\label{compo}
\delta_i\left(\delta_{a}(v)\right)=\left(\delta_{i}^{-}(\delta_{a}(v)),\delta_{a}(v)_{i}\right)=
\left(\delta_{a_i}^-(v),v_{a_i(k+1)}\right)=\delta_{a_i}(v)\in \NN^{m-k}
\end{equation}
 In this notation we have
\begin{equation}\label{compo-2}
 \pi\left(\delta_{a}(v)\right)=\delta_{a}^{-}(v)\quad\mbox{\rm and}\quad
 \pi^+\left(\delta_{a}(v)\right) =a(k)
 \end{equation}
In particular, when $k=1$ for any $a\in \langle m\rangle^{\langle 1\rangle}$, $v\in \NN^{m+1}$  we have
 $$
 \delta_{a}(v)=\left( \delta_{a}^-(v),v_{a(1)}\right)\Longrightarrow
  \pi^+\left(\delta_{a}(v)\right)=v_{a(1)} $$
In addition,  the $i$-th coordinate of
$$
\delta_{a}(v)=(v_1,\ldots,v_{a(1)-1},v_{a(1)+1},\ldots,v_m,v_{a(1)})\in \NN^{m}
$$
with $i\in [m-1]$ is given by some index $b(i)\in [m]-\{a(1)\}$. If we set
$$
a_i(1)=a(1)\quad\mbox{\rm and}\quad a_i(2)=b(i)
$$
we define a mapping $a_i\in \langle m\rangle^{\langle 2\rangle}$ such that
\begin{eqnarray*}
  \delta_{i}(\delta_{a}(v))&=&\left(\delta_{i}^{-}(\delta_{a}(v)),v_{a_i(2)}\right)\\
  &=&
  (v_1,\ldots,v_{a(1)\wedge a_i(2)-1},v_{a(1)\wedge a_i(2)+1},\ldots,v_{a(1)\vee a_i(2)-1},v_{a(1)\vee a_i(2)+1},\ldots,v_m,v_{a_i(2)})\\
  &=&\delta_{a_i}(v)
  \end{eqnarray*}
The above lemma now implies that for any $m>1$ and $v\in\NN^{m+1}$
\begin{eqnarray*}
\displaystyle  M(v)&=&P^{1+\pi^+(v)}~\left[
\tr\left( M( \pi(v))\right)~I~+~2~\sum_{a\in \langle m\rangle^{\langle 1\rangle}}~
 M(\delta_{a}(v))\right]\\
&=&\displaystyle P^{\overline{v}_{m+1}}~\left[
\tr\left( M(\pi(v))\right)~I~+~2~\sum_{a\in \langle m\rangle^{\langle 1\rangle}}~
P^{\overline{v}_{a(1)}}~\tr\left( M(\delta_{a}^{-}(v))\right)\right.\\
&&
\hskip7cm\left.\displaystyle+2^2~\sum_{a\in \langle m\rangle^{\langle 2\rangle} }~P^{\overline{v}_{a(1)}}~
 M(\delta_{a}(v)) \right]
\end{eqnarray*}
Iterating this procedure we find that
\begin{eqnarray*}
\displaystyle  M(v)&=&\displaystyle   ~P^{\overline{v}_{m+1}}~\left[
\tr\left( M(\pi(v))\right)~I~+~\sum_{1\leq k\leq m-2}~2^k~\sum_{a\in \langle m\rangle^{\langle k\rangle}}~
P^{\overline{v}(a)}~\tr\left( M(\delta_{a}^{-}(v))\right)\right]\\
&&
\displaystyle+P^{\overline{v}_{m+1}}~2^{m-1}~\sum_{a\in \langle m\rangle^{\langle m-1\rangle} }~
P^{\overline{v}(a)}
 M(\delta_{a}(v))
\end{eqnarray*}
We end the proof recalling that for any $a\in \langle m\rangle^{\langle m-1\rangle}$ we have
$$
[m]-a[m-1]=\{j\}:=a_j(m)\Longrightarrow
\delta_{a}(v)=(v_{a_j(m)},v_{a_j(m-1)})\quad \mbox{\rm with}\quad a_j(m-1)=a(m-1)
$$
so that
$$
 M(\delta_{a}(v))=  ~P^{1+v_{a(m-1)}}~\tr\left(P^{1+v_{j}}\right)~+~2~  ~P^{1+v_{a_j(m)}}~
P^{1+v_{a_j(m-1)}}
$$
This yields the decomposition
\begin{eqnarray*}
\displaystyle  M(v)&=&\displaystyle   ~P^{\overline{v}_{m+1}}~\left[
\tr\left( M(\pi(v))\right)~I~+~\sum_{1\leq k\leq m-2}~2^k~\sum_{a\in \langle m\rangle^{\langle k\rangle}}~
P^{\overline{v}(a)}~\tr\left( M(\delta_{a}^{-}(v))\right)\right]\\
&&
\displaystyle+  ~P^{\overline{v}_{m+1}}~2^{m-1}~\sum_{a\in \langle m\rangle^{\langle m-1\rangle} }~
P^{\overline{v}(a)}~
\tr\left( M(\delta_{a}^{-}(v))\right)
\displaystyle+P^{\overline{v}_{m+1}}~2^{m}~  
\sum_{a\in \langle m\rangle^{\langle m\rangle} }~
P^{\overline{v}(a)}~
\end{eqnarray*}
This ends the proof of the proposition.
\cqfd

Now we are in position to prove the polynomial formulae (\ref{S-formula}).

{\bf Proof of (\ref{S-formula}) : }\label{S-formula-proof}

 For any $a\in \langle m\rangle^{\langle k\rangle}$
we have
$$
\delta_{a}^{-}(v)=v_b:=(v_{b(1)},\ldots,v_{b(m-k)})
$$
with $b\in \langle m\rangle^{\langle m-k\rangle}$ given by
$$
\left(b(1),\ldots,b(m-k)\right)=\left(1,\ldots,a(\sigma(1))-1,a(\sigma(1))+1,\ldots,a(\sigma(k))-1,a(\sigma(k))+1,\ldots, m\right)
$$
This implies that
$$ 
\overline{v}(a)+\overline{v}(b)=m+\vert v\vert\Longrightarrow
\overline{v}(a)=m+\vert v\vert-\overline{v}(b)=m+\vert v\vert-(m-k)-v(b)=\vert v\vert+k-v(b)
$$
Up to a change of index this yields
$$
\begin{array}{l}
\displaystyle  M(v,l)-2^{m}~m!~~
P^{m+\vert v\vert+l+1}\\
\\
\displaystyle=
\displaystyle ~P^{1+l}~\left[
\tr\left( M(v)\right)~I~+~\sum_{1\leq k<m}~2^{k}~\frac{k!}{(m-k)!}~P^k
\sum_{b\in \langle m\rangle^{\langle m-k\rangle}}~
P^{\vert v\vert-v(b)}~\tr\left( M(v_b)\right)\right]\\
\\
\displaystyle=
~P^{1+l}~\left[
\tr\left( M(v)\right)~I~+~\sum_{1\leq k<m}~2^{m-k}~\frac{(m-k)!}{k!}~P^{m-k}
\sum_{b\in \langle m\rangle^{\langle k\rangle}}~
P^{\vert v\vert-v(b)}~\tr\left( M(v_b)\right)\right]\\
\\
\displaystyle=
~P^{1+l}~\left[
\sum_{1\leq k\leq m}~2^{m-k}~\frac{(m-k)!}{k!}~P^{m-k}
\sum_{b\in \langle m\rangle^{\langle k\rangle}}~
P^{\vert v\vert-v(b)}~\tr\left( M(v_b)\right)\right]\
\end{array}
$$
This implies that
\begin{equation}\label{ref-prop-1}
M(v,l)=
\sum_{0\leq k\leq m}~2^{m-k}~\frac{(m-k)!}{k!}~P^{(m+1)+l+\vert v\vert-k}
\sum_{b\in \langle m\rangle^{\langle k\rangle}}~
P^{-v(b)}~\tr\left( M(v_b)\right)
\end{equation}
with the convention
$$
k=0\Longrightarrow \sum_{b\in \langle m\rangle^{\langle 0\rangle}}~
P^{-v(b)}~\tr\left( M(v_b)\right)=I
$$

In particular when $(v,l)\in\NN^{m+1}$, and $n\in \NN$ are chosen  such that
$
\vert v\vert=n
$
we find the formula
\begin{equation}\label{ref-prop}
M(v,l)=
\sum_{0\leq k\leq m}~2^{m-k}~\frac{(m-k)!}{k!}~P^{1+l+n+m-k}
\sum_{b\in \langle m\rangle^{\langle k\rangle}}~
P^{-v(b)}~\tr\left( M(v_b)\right)
\end{equation}

In addition, for any $v\in \NN^{m+1}$ with $\vert v\vert=n$ we have
\begin{eqnarray*}
M(v)&=&
\sum_{0\leq k\leq m}~2^{k}~(m)_{m-k}~\frac{k!}{(m-k)!}~P^{k+v_{m+1}+1}~~\frac{1}{(m)_{m-k}}
\sum_{b\in \langle m\rangle^{\langle m-k\rangle}}~
P^{n-v(b)}~\tr\left( M(v_b)\right)\\&=&
\sum_{0\leq k\leq m}~2^{k}~(m)_k~P^{k+\pi^+(v)+1}~(m)_{m-k}^{-1}
\sum_{b\in \langle m\rangle^{\langle m-k\rangle}}~
P^{n-v(b)}~\tr\left( M(v_b)\right)\\
&=&\sum_{0\leq q\leq m+n}~\left[\sum_{(k,l)\in \Delta_{m,n}(q)}~2^{k}~(m)_k~t_{\pi(v)}(m-k,n-l)\right]~P^{q+\pi^+(v)+1}
\end{eqnarray*}
This ends the proof of the l.h.s. formula in (\ref{S-formula}).\cqfd

\section{A non commutative binomial formulae}

\subsection{Quadratic forms and extended binomial coefficients}
This section is mainly concerned with the proof of the first part of theorem~\ref{theo-2}.
\begin{defi}
For any $k,l,n\geq 0$, and  $(Q,x)\in(\SS_{r}\times\RR^r)$ we set
$$
Q^{[n]}(x):=\sum_{k+l=n}Q^{[k,l]}(x)\quad\mbox{with the matrix functionals}\quad
Q^{[k,l]}(x):=Q^k\,
{x}\,{x}^{\prime}\, Q^l
$$
In the above display, the summation is taken over integers $k,l\geq 0$.  
We also consider the 
quadratic
polynomial functions defined by the normalized trace formula
\begin{equation}\label{def-s}
s_n(x):=\frac{1}{n+1}~\tr\left(Q^{[n]}(x)\right)=\langle x,Q^{n}x\rangle
\end{equation}
\end{defi}
Next lemma provides an recursive formula for computing sequentially 
the quadratic form functionals $Q^{[n]}$.

\begin{lem}
We have $Q^{[0]}(x)=\,xx^{\prime}$ and for any $n\geq 0$
\begin{equation}\label{rec-PP-x}
2Q^{[n+1]}(x)-\left[Q^{[n+1,0]}(x)+Q^{[0,n+1]}(x)\right]=
Q~Q^{[n]}(x)+Q^{[n]}(x)~Q
\end{equation}

\end{lem}
\proof
We have
$$
\sum_{k+l=n}\left[Q^{[k+1,l]}(x)+Q^{[k,l+1]}(x)\right]=
Q~Q^{[n]}(x)+Q^{[n]}(x)~Q
$$
and for any $n\geq 0$ we have
\begin{eqnarray*}
2Q^{[n+1]}(x)&=&\sum_{k+l=n}\left[Q^{k+1}\,xx^{\prime}\,Q^{l}+P^{k}\,xx^{\prime}\,Q^{l+1}\right]+Q^{n+1}\,xx^{\prime}\,+\,xx^{\prime}\,Q^{n+1}\\
&=&Q~Q^{[n]}(x)+Q^{[n]}(x)~Q+
Q^{[n+1,0]}(x)+Q^{[0,n+1]}(x)
\end{eqnarray*}
This ends the proof of the lemma.
\cqfd

\begin{defi}
For any given parameter $n$ let 
$$
\forall 0\leq k\leq n\qquad \Xi_{k,n}~:~(l,x)\in\left(\NN\times \RR^r\right)~\mapsto ~\Xi_{k,n}^l(x)\in~[0,\infty[
$$
be the functions with support $\Ia_{k,n}=\left(\{0,\ldots,k\}\times \RR^r\right)$ defined for any $k<n$ by the
 the backward recursions 
$$
\forall  (l,x)\in \left(\NN\times \RR^r\right)\qquad\begin{array}[t]{rcl}
 \Xi_{k,n}^l(x)
 &=&\displaystyle\sum_{l<j}~ \Xi_{k+1,n}^j(x)~s_{j-(l+1)}(x)
 \end{array}
$$
with the constant boundary condition function $ \Xi^l_{n,n}(x)=1_{\{n\}}(l)$. 
 \end{defi}
 In terms of extended binomial coefficients we have
\begin{eqnarray}
 \Xi_{k,n}^l(x)&=& \Xi_{l+(k-l),n}^l(x)=\Xi_{k-l,n-l}^0(x)
 \nonumber\\
&& \nonumber\\
 &=&\sum_{v_1+\ldots+v_{n-k}=k-l}~s_{v_1}(x)~\ldots s_{v_{n-k}}(x):=\left(
\begin{array}{c}
(n-l)-1\\
k-l
\end{array}
\right)_{s(x)}\label{generalized-binomial}
\end{eqnarray}

with the conventions
$$
\left(
\begin{array}{c}
-1\\
k
\end{array}
\right)_{s(x)}=\Xi^0_{k,0}=1_{k=0}=\Xi^0_{k,k}=\left(
\begin{array}{c}
k-1\\
k
\end{array}
\right)
\quad\mbox{and}\quad
\left(
\begin{array}{c}
n-1\\
0
\end{array}
\right)_{s(x)}= \Xi_{0,n}^0=s_0(x)^{n}
$$

The equivalent formula (\ref{generalized-binomial}) comes from the fact that
$$
\begin{array}[t]{rcl}
 \Xi_{k,n}^l
 &=&\displaystyle\sum_{j_1\geq 0}~s_{j_1}~ \Xi_{k+1,n}^{(l+1)+j_1}=\displaystyle\sum_{j_1,j_2\geq 0}~s_{j_1}~s_{j_2}~ \Xi_{k+2,n}^{(l+2)+j_1+j_2} \end{array}
$$
Iterating this procedure we check  (\ref{generalized-binomial}).

The above formula shows that $ \Xi_{k,n}^l$ only depends on the parameters $(n-k)$ and $(k-l)$.

 \begin{lem}
 For any $n\geq 0$ and $l\leq n$ we have 
 $
   \Xi_{n,n+1}^{l}=s_{n-l}
 $.
  In addition, for any $k,l\in \NN$ we have
 \begin{equation}\label{ref-2-binomial}
 \Xi_{k+1,n+1}^{l+1}= \Xi_{k,n}^l\qquad\mbox{and}\quad k<n~\Longrightarrow~
  \Xi_{k,n}^k(x)=s_{0}(x)^{n-k}
 \end{equation}
 
 \end{lem}
 \proof

 Observe that for any $k<n$ we have
 $$
  \Xi_{k,n}^k(x)
 =\displaystyle\sum_{k<j}~ \Xi_{k+1,n}^j(x)~s_{j-(k+1)}(x)= \Xi_{k+1,n}^{k+1}(x)~s_{0}(x)\Longrightarrow  \Xi_{k,n}^k=s_{0}^{n-k}
  $$

For any $n\geq 1$ we also have
 $$
  \Xi_{n-1,n}^0(x)=\displaystyle\sum_{0<j}~ \Xi_{n,n}^j(x)~s_{j-1}(x)=s_{n-1}(x)
 $$
 In addition, for any $1\leq l\leq n$ we have
 $$
  \Xi_{n,n+1}^{l}(x)
 =\displaystyle\sum_{l<j}~ \Xi_{n+1,n+1}^j(x)~s_{j-(l+1)}(x)=s_{n-l}(x)=
 s_{(n-1)-(l-1)}(x)=  \Xi_{n-1,n}^{l-1}(x)
  $$
This also shows that $\Xi_{n,n+1}^{0}=s_{n}$. More generally for any $k,l\in\NN$ we have
$$
 \Xi_{k+1,n+1}^{l+1}= \Xi_{k,n}^l
 $$
We prove this assertion using a backward induction w.r.t. the parameter $k$.
The result has been checked for $k=n-1$. Assume it has been checked at rank $k+1$. In this case we have
\begin{eqnarray*}
 \Xi_{k,n}^{l+1}(x)
 &=&\displaystyle\sum_{l+1<j}~ \Xi_{k+1,n}^j(x)~s_{j-(l+2)}(x)\\
 &=&
  \displaystyle\sum_{l< j}~ \Xi_{k+1,n}^{j+1}(x)~s_{j-(l+1)}(x)=
   \displaystyle\sum_{l< j}~  \Xi_{k,n}^j(x)~s_{j-(l+1)}(x)=\Xi_{k-1,n}^l(x)
\end{eqnarray*}
This ends the proof of the lemma.
\cqfd

\subsection{An almost sure binomial formula}
This section is mainly concerned with the proof of the following theorem.
\begin{theo}
For any $(n,x,Q)\in(\NN\times\RR^r\times\SS_r)$ we have the non commutative binomial
formula
\begin{eqnarray}
\left(xx^{\prime}+Q\right)^{n+1}
&=&Q^{n+1}+
\sum_{0\leq l\leq k\leq n}~
\left(
\begin{array}{c}
(n-l)-1\\
k-l
\end{array}
\right)_{s(x)}~
Q^{[l]}(x)~\label{explicit-xx-tP}
\end{eqnarray}
with the function $s(x)$ and binomial coefficients defined in (\ref{def-s}) and (\ref{generalized-binomial}).
\end{theo}
\proof
We prove (\ref{explicit-xx-tP}) by induction w.r.t. the parameter $n$. For $n=0$
the proof of the formula is immediate. We assume that it has been proved at some rank
$n\geq 0$. We use the decomposition
$$
\begin{array}{l}
\displaystyle2\left[\left(xx^{\prime}+Q\right)^{n+2}-~Q^{n+2}\right]\\
\\
\displaystyle=\left(xx^{\prime}+Q\right)^{n+1}\left(xx^{\prime}+Q\right)+
\left(xx^{\prime}+Q\right)\left(xx^{\prime}+Q\right)^{n+1}-2Q^{n+2}\\
\\
\displaystyle=\left[xx^{\prime}Q^{n+1}+Q^{n+1}xx^{\prime}\right]+
\sum_{0\leq k\leq n}~\sum_{l=0}^k~
\Xi_{k,n}^l(x)~
\left[xx^{\prime}Q^{[l]}(x)+Q^{[l]}(x)xx^{\prime}\right]\\
\\
\displaystyle\hskip3cm+\sum_{0\leq k\leq n}~\sum_{l=0}^k~
\Xi_{k,n}^l(x)~
\left[QQ^{[l]}(x)+Q^{[l]}(x)Q\right]
\end{array}
$$
Observe that
\begin{eqnarray*}
xx^{\prime}Q^{[l]}(x)+Q^{[l]}(x)xx^{\prime}&=&\sum_{u+v=l}xx^{\prime}Q^u
{x}\,{x}^{\prime}\, Q^v\\
&=&\sum_{u+v=l}~ s_u(x)  \left[Q^{[v,0]}(x)+Q^{[0,v]}(x)\right]
\end{eqnarray*}
Using (\ref{rec-PP-x}) this implies that
$$
\begin{array}{l}
\displaystyle2\left[\left(xx^{\prime}+Q\right)^{n+2}-Q^{n+2}\right]\\
\\
\displaystyle=\left[Q^{[n+1,0]}(x)+Q^{[0,n+1]}(x)\right]+
\sum_{0\leq k\leq n}\sum_{v=0}^k\left[
\sum_{l=v}^k~\Xi_{k,n}^l(x)~ s_{l-v}(x) \right]~ \left[Q^{[v,0]}(x)+Q^{[0,v]}(x)\right]\\
\\
\displaystyle\hskip3cm+\sum_{0\leq k\leq n}~t^{k+1}~\sum_{l=0}^k~
\Xi_{k,n}^l(x)~
\left[2Q^{[l+1]}(x)-\left[Q^{[l+1,0]}(x)+Q^{[0,l+1]}(x)\right]\right]
\end{array}
$$
Simplifying the first term and applying twice (\ref{ref-2-binomial}) we have
$$
\begin{array}{l}
\displaystyle2\left[\left(xx^{\prime}+Q\right)^{n+2}-Q^{n+2}\right]\\
\\
\displaystyle=~2Q^{[n+1]}(x)+
\sum_{0\leq k\leq n}~\sum_{l=0}^k\left[
\sum_{l<j}
 \Xi_{k+1,n+1}^{j} (x)  ~ 
s_{j-(l+1)}(x) \right]~ \left[Q^{[l,0]}(x)+Q^{[0,l]}(x)\right]\\
\\
\displaystyle\hskip3cm+\sum_{1\leq k\leq  n}~\sum_{l=1}^{k}~ \Xi_{k,n+1}^{l}
\left[2Q^{[l]}(x)-\left[Q^{[l,0]}(x)+Q^{[0,l]}(x)\right]\right]
\end{array}
$$
using the backward induction formula, this yields the formula
$$
\begin{array}{l}
\displaystyle2\left[\left(xx^{\prime}+Q\right)^{n+2}-Q^{n+2}\right]\\
\\
\displaystyle=2Q^{[n+1]}(x)+2~\Xi_{0,n+1}^{0}(x)~Q^{[0]}(x)
+\sum_{1\leq k\leq n}~\sum_{l=0}^k~\Xi_{k,n+1}^{l}(x)~ \left[Q^{[l,0]}(x)+Q^{[0,l]}(x)\right]\\
\\
\displaystyle\hskip3cm+\sum_{1\leq k\leq  n}~~\sum_{l=1}^{k}~ \Xi_{k,n+1}^{l}(x)
\left[2Q^{[l]}(x)-\left[Q^{[l,0]}(x)+Q^{[0,l]}(x)\right]\right]\\
\\
\displaystyle=2Q^{[n+1]}(x)+2~\Xi_{0,n+1}^{0}(x)~Q^{[0]}(x)+2\sum_{1\leq k\leq  n}~~\sum_{l=1}^{k}~ \Xi_{k,n+1}^{l}(x)
Q^{[l]}(x)
\end{array}
$$
This ends the proof of the theorem.
\cqfd

\subsection{Central matrix moments}\label{proof-sec-ae-Binomial-W}

This section is mainly concerned with the proof of the non commutative binomial formulae
(\ref{ae-Binomial-W}) and (\ref{Binomial-W}). 
stated in the first part of theorem~\ref{theo-2}.

Using (\ref{explicit-xx-tP}) for any $n\geq 1$ we have
$$
\begin{array}{l}
\displaystyle\left(\Xa-P\right)^{n}\\
\\
\displaystyle=(-1)^{n}~P^{n}+\sum_{0\leq k< n}\sum_{0\leq 
l\leq k}
\left(
\begin{array}{c}
(n-1)-(l+1)\\
k-l
\end{array}
\right)_{s(X)}~(-1)^{k-l}~\times~(-1)^{l}~
P^{[l]}(X)\\
\\
\displaystyle=(-1)^{n}~P^{n}+\sum_{0\leq k< n}~(-1)^{k}~
\sum_{ 0\leq l\leq k}
\left(
\begin{array}{c}
(n-1)-(l+1)\\
k-l
\end{array}
\right)_{s(X)}
P^{[l]}(X)
\end{array}
$$
with the trace mapping
$
s_n(x):=\langle x,P^{n}x\rangle$. This ends the proof of (\ref{ae-Binomial-W}).

Taking the expectation, we also have
$$
\begin{array}{l}
\displaystyle
\EE\left[\left(
\begin{array}{c}
(n-1)-(l+1)\\
k-l
\end{array}
\right)_{s(X)}
P^{[l]}(X)\right]\\
\\
\displaystyle=~\sum_{u+w=l}~
\EE\left[\sum_{v_1+\ldots+v_{(n-1)-k}=k-l}
P^{u}\,\Xa \,P^{v_1}\,\Xa \ldots P^{v_{(n-1)-k}}
\,\Xa \, P^{w}\right]\\
\\
\displaystyle=~\sum_{u+w=l}~\sum_{v_1+\ldots+v_{(n-1)-k}=k-l}
P^{u}~M(v_1,\ldots,v_{(n-1)-k},w)=\sum_{u+w=l}~\sum_{v\in V_{(n-1)-k,k-l}}
P^{u}~M(v,w)
\end{array}
$$
Observe that
$$
k=l=0=u=w\Rightarrow \sum_{v\in V_{(n-1)-k,k-l}}
P^{u}~M(v,w)
=\EE(\Xa^n)
$$
By (\ref{f1}) for any $w+u=l$ we also have
\begin{eqnarray*}
P^{u}~M(v,w)&=&M(v,w)~P^{u}=M(v,w+u)=
M(v,l)
\end{eqnarray*}
This yields
$$
\begin{array}{l}
\displaystyle\EE\left[\left(\Xa-P\right)^{n}-\Xa^n\right]-(-1)^{n}~P^{n}\\
\\
\displaystyle=\sum_{1\leq k<n}~(-1)^k~\sum_{0\leq 
l\leq k}~\sum_{v_1+\ldots+v_{(n-1)-k}+v_{n-k}=k}~~(1+v_{n-k})~
M(v_1,\ldots,v_{n-k})\\
\\
\displaystyle=\sum_{1\leq k< n}~(-1)^k~W(n-k,k)
\end{array}
$$
This ends the proof of (\ref{Binomial-W}).
When $n=2m$ we have
$$
\begin{array}{l}
\displaystyle\EE\left[\left(\Xa-P\right)^{2m}-\Xa^{2m}\right]\\
\\
\displaystyle=P^{2m}+\sum_{1\leq k\leq m-1}~W(2(m-k),2k)-
\sum_{1\leq k\leq m}~W(2(m-k)+1,2k-1)\\
\\
\displaystyle=-(2m-1)~P^{2m}-\sum_{1\leq k<m}\left[W(2(m-k)+1,2k-1)-W(2(m-k),2k)\right]
\end{array}
$$
The last assertion comes from the fact that
$$
W(1,2m-1)=2m~M(2m-1)=2m~P^{2m}
$$
When $n=2m+1$ we have
$$
\begin{array}{l}
\displaystyle\EE\left[\left(\Xa-P\right)^{2m+1}-\Xa^{2m+1}\right]+~P^{2m+1}\\
\\
\displaystyle=\sum_{1\leq k\leq 2m}~(-1)^k~W(2m+1-k,k)\\
\\
\displaystyle=-\sum_{1\leq k\leq m}~\left[W(2(m-k)+2,2k-1)-W(2(m-k)+1,2k)\right]
\end{array}
$$

This ends the proof of the (\ref{ae-Binomial-W}) and (\ref{Binomial-W}).
\cqfd
\subsection{Weighted matrix product moments}\label{proof-sec-weighted}

This section is mainly concerned with the proof of the r.h.s. polynomial formulae
stated in (\ref{S-formula}).

By (\ref{ref-prop}) we find that
$$
\begin{array}{l}
\displaystyle  \sum_{v_1+\ldots+v_m+l=n}~(l+1)~M(v,l)\\
\\
=\displaystyle \sum_{0\leq l\leq n}~(l+1)~\sum_{v_1+\ldots+v_m= {n-l}}~M(v,l)\\
\\
=\displaystyle\sum_{0\leq k\leq m}\sum_{0\leq l\leq n}~(l+1)~~2^{m-k}~\frac{(m-k)!}{k!}~P^{1+l+ {n-l}+m-k}\\
\\
\displaystyle\hskip3cm~\sum_{b\in \langle m\rangle^{\langle k\rangle}}~\sum_{v_{b(1)}+\ldots+v_{b(k)}+\sum_{k<q\leq m}v_q=n-l}~
P^{l-n+\sum_{k<q\leq m}v_q}~\tr\left( M(v_b)\right)
\end{array}
$$
This implies that
$$
\begin{array}{l}
\displaystyle  \sum_{v_1+\ldots+v_m+l=n}~(l+1)~M(v,l)\\
\\
=\displaystyle\sum_{0\leq k\leq m}\sum_{0\leq l_1+l_2\leq n}~(l_1+1)~~2^{m-k}~\frac{(m-k)!}{k!}~P^{1+l_1+l_2+m-k}\\
\\
\displaystyle\hskip3cm~(m)_k~\left[\sum_{v_{k+1}+\ldots+v_m=l_2}~1\right]
\sum_{v_{1}+\ldots+v_{k}=n-l_1-l_2}~
\tr\left( M(v_1,\ldots,v_k)\right)\\
\\
=\displaystyle\sum_{0\leq k\leq m}\sum_{0\leq l_1+l_2\leq n}~(l_1+1)~~2^{m-k}~\frac{(m-k)!}{k!}~P^{1+l_1+l_2+m-k}\\
\\
\displaystyle\hskip3cm~(m)_k~\left(
\begin{array}{c}
(m-k)+l_2-1\\
l_2
\end{array}
\right)~\tau(k,n-(l_1+l_2))
\end{array}
$$
Changing the index we conclude that
$$
\begin{array}{l}
\displaystyle W(m+1,n)
=\displaystyle\sum_{0\leq k\leq m}\sum_{0\leq l_1+l_2\leq n}~(l_1+1)~~2^{k}~\frac{m!}{(m-k)!}~P^{1+l_1+l_2+k}\\
\\
\displaystyle\hskip3cm~\left(
\begin{array}{c}
k+l_2-1\\
l_2
\end{array}
\right)~\tau(m-k,n-(l_1+l_2))
\end{array}
$$
In summary we have
$$
\begin{array}{l}
\displaystyle W(m+1,n)
=\displaystyle\sum_{0\leq k\leq m}\sum_{0\leq l\leq n}
~2^{k}~\frac{m!}{(m-k)!}~P^{1+l+k}\\
\\
\displaystyle\hskip3cm~\left[\sum_{l_1+l_2=l}~(l_1+1)\left(
\begin{array}{c}
k+l_2-1\\
l_2
\end{array}
\right)\right]~\tau(m-k,n-l)
\end{array}
$$
 Recalling that
 $$
 \sum_{0\leq k\leq n}\left(
\begin{array}{c}
q+k\\
q
\end{array}
\right)=\left(
\begin{array}{c}
q+n+1\\
n
\end{array}
\right)\Longrightarrow
 \sum_{0\leq k\leq n}~(n-k)~\left(
\begin{array}{c}
q+k\\
q
\end{array}
\right)=\left(
\begin{array}{c}
q+n+1\\
n-1
\end{array}
\right)
 $$
 We check the r.h.s. assertion using the fact that
 $$
~k~\left(
\begin{array}{c}
q+k\\
q
\end{array}
\right)=(q+1)~\left(
\begin{array}{c}
q+1+(k-1)\\
q+1
\end{array}
\right)
 $$
This implies that
\begin{eqnarray*}
\sum_{l_1+l_2=l}~(l_1+1)\left(
\begin{array}{c}
k+l_2-1\\
l_2
\end{array}
\right)&=&\sum_{0\leq l_2\leq l+1}~((l+1)-l_2)\left(
\begin{array}{c}
(k-1)+l_2\\
(k-1)
\end{array}
\right)\\
&=&\left(
\begin{array}{c}
(k-1)+(l+2)\\
l
\end{array}
\right)=\left(
\begin{array}{c}
k+l+1\\
l
\end{array}
\right)
\end{eqnarray*}
This implies that
$$
\begin{array}{l}
\displaystyle W(m+1,n)\\
\\
=\displaystyle\sum_{0\leq k\leq m}\sum_{0\leq l\leq n}
~2^{k}~\frac{m!}{(m-k)!}~\left(
\begin{array}{c}
k+l+1\\
l
\end{array}
\right)\tau(m-k,n-l)~~P^{1+k+l}\\
\\
=\displaystyle\sum_{0\leq q\leq m+n}~ \left[\sum_{(k,l)\in \Delta_{m,n}(q)}
~2^{k}~(m)_k~~~\left(
\begin{array}{c}
k+l+1\\
l
\end{array}
\right)
\tau(m-k,n-l)\right]~~P^{1+q}
\end{array}
$$
This ends the proof of the l.h.s.  formulae
stated in (\ref{S-formula}).\cqfd
\section{Some moments estimates}\label{sec-m-e}
\begin{defi}\label{def-theta-i}
For any $1\leq i\leq m$ we let $\theta_i$ and $\theta_i^+:v\in\NN^{m+1}\mapsto \theta_i^+(v)\in\NN^m$  be the contraction mappings
defined by
\begin{eqnarray*}
\theta_i(v_1,\ldots,v_{m+1})&=&\left(v_1,\ldots,v_{i-1},1+v_{m+1}+v_{i},v_{i+1},\ldots,,v_m\right)\\
\theta_i^+(v_1,\ldots,v_{m+1})&=&\left(v_1,\ldots,v_{i-1},v_{i+1},\ldots,,v_m,1+v_{m+1}+v_{i}\right)
\end{eqnarray*}

\end{defi}

In this notation, for any $v\in \NN^{m+1}$ and $1\leq i\leq m$ we have
\begin{equation}\label{delta-theta}
P^{1+\pi^+(v)}~M(\delta_i(v))=M(\theta_i^+(v))
\end{equation}
By the symmetry of the trace 
$$
\forall \sigma\in\Ga_m\qquad \tr\left(M(v_1,\ldots,v_m)\right)=\tr\left(M(v_{\sigma(1)},\ldots,v_{\sigma(m)})\right)
$$
we also have
\begin{equation}\label{delta-theta-bis}
\sum_{1\leq i\leq m}\tr\left(M(\theta_i^+(v))\right)=\sum_{1\leq i\leq m}\tr\left(M(\theta_i(v))\right)
\end{equation}

\begin{lem}
For any $v\in \NN^{m+1}$ we have
\begin{equation}\label{lower-bound}
\left(1+\frac{1}{2m}~\right)~\sum_{1\leq i\leq m}
M\left(\theta_i^+(v)\right)\leq ~\frac{1}{2}~
 M(v)
\end{equation}

\end{lem}
\proof

By (\ref{delta-theta}) he recursion (\ref{recursion-XX-moments-intro}) implies that
\begin{eqnarray*}
\displaystyle 
M(v)
&=&P^{1+\pi^+(v)}~\left[\tr\left(M(\pi(v))\right)~I+2~\sum_{1\leq i\leq m}~M(
\delta_i(v))\right]\\
&=& \tr\left(M(\pi(v))\right)~P^{1+\pi^+(v)}+2~\sum_{1\leq i\leq m}~M(\theta_i^+(v))\geq M(\theta_j^+(v))+2~\sum_{1\leq i\leq m}~M(\theta_i^+(v))
\end{eqnarray*}
The last assertion comes from the  symmetry of the trace and the fact that $\tr\left(M(\pi(v))\right)~I\geq M(\pi(v)$.
Summing over all $1\leq j\leq m$ we find that
$$
\forall v\in \NN^{m+1}\qquad
\displaystyle ~ M(v)
 \geq \displaystyle ~2~\left(1+\frac{1}{2m}\right)~~\sum_{1\leq j\leq m}
 M\left(\theta_j^+(v)\right)
$$
The proof of the lemma is now completed.
\cqfd

\begin{lem}\label{F-lemma}
Let $F$ be some function from $\cup_{m\geq 1}\NN^m$ into $[0,\infty[$ 
Assume that  for any $v\in \NN^{m+1}$ we have
\begin{equation}\label{Hyp-F}
F(v) \geq 2~\left(1+\epsilon_m\right)\sum_{1\leq j\leq m}~F\left(\theta_j(v)\right)
\end{equation}
for some $\epsilon_m\in [0,1]$.
In this situation we have
\begin{equation}\label{first-estimate}
\begin{array}{l}
\displaystyle 
\sum_{v\in V_{m,n+1}}~F(v)
\displaystyle \leq  \frac{1}{2}~\left(1-\frac{\epsilon_m}{1+\epsilon_m}\right) \sum_{v\in V_{m+1,n}}~F(v)
 \end{array}
\end{equation}
as well as
\begin{equation}\label{Hyp-F-equation}
\left(1+\epsilon_m\right)~\sum_{v\in V_{m,n+1}}~(1+\pi^+(v))~F(v)\leq \sum_{v\in V_{m+1,n}}~(1+\pi^+(v))~F(v)
\end{equation}
\end{lem}

The proof of the above technical lemma is housed in the appendix on page~\pageref{proof-F-lemma}.

Next proposition is pivotal in the proof of the estimates (\ref{rhoW-inq}) stated in theorem~\ref{theo-1}.

\begin{prop}\label{prop-tau-varpi}
For any $m\geq 1$ and any $n\geq 0$ we have
$$
\displaystyle 2~\left(1+\frac{1}{m}\right) \tau(m,n+1)
\displaystyle \leq ~\tau(m+1,n)
\quad
and
\quad
\left(1+\frac{1}{m}\right)~\varpi(m,n+1)\leq \varpi(m+1,n)
$$
\end{prop}

\proof
Taking the trace in (\ref{lower-bound}) and using (\ref{delta-theta-bis}) we have
$$
\forall v\in \NN^{m+1}\qquad
\displaystyle 
\tr\left(M(v)\right)
\displaystyle \geq 2~\left(1+\frac{1}{2m}\right)~\sum_{1\leq j\leq m}~
 \tr\left[M\left(\theta_j(v)\right)\right]
$$
The end  of the proof is now a direct consequence of the lemma~\ref{F-lemma}. The proof of the proposition
is completed
\cqfd

We are now in position to prove the estimates (\ref{rhoW-inq}) stated in the second part of theorem~\ref{theo-1}

{\bf Proof of (\ref{rhoW-inq}):}\label{proof-rhoW-inq}

Observe that
$$
\forall (k,l)\in \Delta_{m-1,n}(q)\qquad q=k+l<m+n
$$
Thus, by (\ref{convention-rho}),  for any $q\leq m+n$ we have
\begin{eqnarray*}
\rho^W_{m,n}(q)&=&\displaystyle\sum_{(k,l)\in \Delta_{m-1,n}(q)}
~2^{k}~(m)_k~~~\left(
\begin{array}{c}
q+1\\
l
\end{array}
\right)
\tau(m-k,n-l)~\\
&&\hskip7cm\displaystyle+~~1_{q=m+n}~~2^{m}~m!~~~\left(
\begin{array}{c}
m+n+1\\
n
\end{array}
\right)
\end{eqnarray*}
By (\ref{defi-rho}) we also have the decomposition
\begin{eqnarray*}
\rho^W_{m-1,n+1}(q)
&=&\displaystyle~1_{q<m+n}~\sum_{(k,l)\in \Delta_{m-1,n}(q)}
~2^{k}~(m-1)_k~~~\left(
\begin{array}{c}
q+1\\
l
\end{array}
\right)
\tau(m-1-k,n+1-l)\\
&&\hskip1.5cm+~1_{q>n}~2^{q-(n+1)}~(m-1)_{q-(n+1)}~~~\left(
\begin{array}{c}
q+1\\
n+1
\end{array}
\right)
\tau(m+n-q,0)
\end{eqnarray*}
This yields the formula
\begin{eqnarray*}
\rho^W_{m-1,n+1}(q)
&=&\displaystyle 1_{q<m+n}~\sum_{(k,l)\in \Delta_{m-1,n}(q)}
~2^{k}~(m-1)_k~~~\left(
\begin{array}{c}
q+1\\
l
\end{array}
\right)
\tau(m-1-k,n+1-l)\\
&&\hskip1.5cm+~1_{n<q<n+m}~2^{q-(n+1)}~(m-1)_{q-(n+1)}~~~\left(
\begin{array}{c}
q+1\\
n+1
\end{array}
\right)
\tau(m+n-q,0)\\
\\
&&\hskip2cm+~1_{q=n+m}~2^{m-1}~(m-1)!~~~\left(
\begin{array}{c}
n+m+1\\
n+1
\end{array}
\right)
\end{eqnarray*}
Comparing the summands associated with $l=n\Longrightarrow k=q-n$ this implies that
$$
\rho^W_{m,n}(q)-\rho^W_{m-1,n+1}(q)=A+B+C
$$
with
$$
\begin{array}{l}
\displaystyle A =1_{q<m+n}~
\sum_{(k,l)\in \Delta_{m-1,n-1}(q)}~2^k~\left(
\begin{array}{c}
q+1\\
l
\end{array}
\right)\\
\\
\hskip3cm\displaystyle\times\left[
~(m)_k~~
\tau(m-k,n-l)-~(m-1)_k~
\tau(m-1-k,n+1-l)\right]~\\
\\
\displaystyle B=1_{q=m+n}~\left[2^{m}~m!~~~\left(
\begin{array}{c}
m+n+1\\
n
\end{array}
\right)-2^{m-1}~(m-1)!~~~\left(
\begin{array}{c}
m+n+1\\
n+1
\end{array}
\right)
\right]
\end{array}
$$
and
$$
\begin{array}{l}
\displaystyle
C=1_{n<q<n+m}\left\{2^{q-n}~\left(
\begin{array}{c}
q+1\\
n
\end{array}
\right)\left[
~(m)_{q-n}~~
\tau(m+n-q,0)-~(m-1)_{q-n}~
\tau(m+n-1-q,1)\right]\right.\\
\\\left.\hskip6cm-~2^{q-(n+1)}~(m-1)_{q-(n+1)}~~~\left(
\begin{array}{c}
q+1\\
n+1
\end{array}
\right)
\tau(m+n-q,0)\right\}
\end{array}
$$
Using proposition~\ref{prop-tau-varpi} we have
$$
\forall n<q<n+m\qquad
\tau(m+n-(q+1),1)\leq \frac{1}{2}~\left(1-\frac{1}{m+n-q}\right)~\tau(m+n-q,0)
$$
and
$$
\forall (k,l)\in \Delta_{m-1,n-1}(q)\qquad
\tau(m-1-k,n+1-l)\leq \frac{1}{2}~\left(1-\frac{1}{m-k}\right)~\tau(m-k,n-l)
$$
The latter yields the estimate
$$
\begin{array}{l}
\displaystyle ~(m)_k~~
\tau(m-k,n-l)-~(m-1)_k~
\tau(m-1-k,n+1-l)\\
\\
\displaystyle\geq  ~(m)_k~\tau(m-k,n-l)~\left[1-\frac{1}{2}~
\left(1-\frac{1}{m-k}\right)~\left(1-\frac{k}{m}\right)\right]\geq  \frac{1}{2}~(m)_k~\tau(m-k,n-l)~
\end{array}
$$
from which we prove that
$$
\begin{array}{l}
\displaystyle A \geq 1_{q<m+n}~\frac{1}{2}~
\sum_{(k,l)\in \Delta_{m-1,n-1}(q)}~2^k~\left(
\begin{array}{c}
q+1\\
l
\end{array}
\right)~(m)_k~\tau(m-k,n-l)\end{array}
$$
In the same vein, we have
$$
\begin{array}{l}
\displaystyle C\geq 1_{n<q<n+m}~2^{q-n}~\tau(m+n-q,0)~\left(
\begin{array}{c}
q+1\\
n
\end{array}
\right)\\
\\
\displaystyle\times\left\{
~(m)_{q-n}~~
-~(m-1)_{q-n}~
\frac{1}{2}~\left(1-\frac{1}{m+n-q}\right)-~\displaystyle \frac{1}{2}~(m-1)_{q-(n+1)}~~\frac{q+1-n}{n+1}\right\}\\
\\
~~=\displaystyle  1_{n<q<n+m}~2^{q-n}~(m)_{q-n}~\tau(m+n-q,0)~\left(
\begin{array}{c}
q+1\\
n
\end{array}
\right)\\
\\
\displaystyle\times\frac{1}{2}~\left\{
\left[1
-\left(1-\frac{q-n}{m}\right)~\left(1-\frac{1}{m+n-q}\right)\right]+\left[1-\displaystyle \frac{1}{n+1}~\left(\frac{1}{m}+\frac{q-n}{m}\right)\right]\right\}\\
\\
\end{array}
$$
This yields the lower bound
$$
\begin{array}{l}
\displaystyle C\geq \displaystyle  1_{n<q<n+m}~\frac{1}{2}~2^{q-n}~~(m)_{q-n}~~\tau(m+n-q,0)~\left(
\begin{array}{c}
q+1\\
n
\end{array}
\right)~\left[1-\displaystyle \frac{1}{n+1}\right]\\
\\
\end{array}
$$
We conclude that for any $m\geq 1$ and $n\geq 0$ we have
$$
\begin{array}{l}
\rho^W_{m,n}(q)-\rho^W_{m-1,n+1}(q)\\
\\
\displaystyle \geq 1_{q<m+n}~\frac{1}{2}~
\sum_{(k,l)\in \Delta_{m-1,n-1}(q)}~2^{k}~\left(
\begin{array}{c}
q+1\\
l
\end{array}
\right)~~(m)_k~\tau(m-k,n-l)~\\
\\
\hskip1cm+\displaystyle~1_{n<q<n+m}~\frac{1}{2}~2^{q-n}~~(m)_{q-n}~~\tau(m+n-q,0)~\left(
\begin{array}{c}
q+1\\
n
\end{array}
\right)~\left[1-\displaystyle \frac{1}{n+1}\right]\\
\\
\hskip4cm\displaystyle+1_{q=m+n}~\frac{1}{2}~2^{m}~m!~\left(
\begin{array}{c}
m+n+1\\
n
\end{array}
\right)~\left[2-\frac{1}{n+1}~\left(1+\frac{1}{m}\right)\right]
\end{array}
$$
This yields  for any $m\geq 1$ and $n\geq 0$ the rather crude estimate
$$
\begin{array}{l}
\rho^W_{m,n}(q)-\rho^W_{m-1,n+1}(q)\\
\\
\displaystyle  \geq \frac{1}{2}~\left[1-\displaystyle \frac{1}{n+1}\right]~
\sum_{(k,l)\in \Delta_{m-1,n}(q)}~2^{k}~\left(
\begin{array}{c}
q+1\\
l
\end{array}
\right)~~(m)_k~\tau(m-k,n-l)~\\
\\
\hskip6cm\displaystyle+1_{q=m+n}~\left[1-\frac{1}{n+1}\right]~2^{m}~m!~\left(
\begin{array}{c}
m+n+1\\
n
\end{array}
\right)~\\
\\
\displaystyle  \geq \frac{1}{2}~\left[1-\displaystyle \frac{1}{n+1}\right]~\rho^W_{m,n}(q)
\end{array}
$$
By (\ref{S-formula}) we also have
\begin{eqnarray*}
 W(m+1,n)- W(m,n+1)&=&\sum_{0\leq q\leq m+n}~\left(\rho^W_{m,n}(q)-\rho^W_{m-1,n+1}(q)\right)~P^{1+q}\\
&\geq&~\frac{1}{2}~\left[1-\displaystyle \frac{1}{n+1}\right]~W(m+1,n)
  \end{eqnarray*}
This ends the proof of the estimates (\ref{rhoW-inq}).
\cqfd

We end this section with the proof of the estimates (\ref{estimate-XX-n}).\label{proof-estimate-XX-n}
By  (\ref{rhoW-inq}) we have $$
\begin{array}{l}
\displaystyle\EE\left[\left(\Xa-P\right)^{2m}-\Xa^{2m}\right]\\
\\
\displaystyle=P^{2m}+\sum_{1\leq k\leq m-1}~W(2(m-k),2k)-
\sum_{1\leq k\leq m}~W(2(m-k)+1,2k-1)\\
\\
\displaystyle\leq -(2m-1)~P^{2m}-
 \sum_{1\leq k<m}\left[W(2(m-k)+1,2k-1)-W(2(m-k),2k)\right]\\
\\
\displaystyle\leq -(2m-1)~P^{2m}-\frac{1}{2}
 \sum_{1\leq k<m}~\left(1-\frac{1}{2k}\right)~W(2(m-k)+1,2k-1)\\
 \\
 \displaystyle\leq -(2m-1)~P^{2m}-\frac{1}{2}~\left(1-\frac{1}{2m}\right)
 \sum_{1\leq k<m}~W(2(m-k)+1,2k-1)
 \end{array}
$$
as soon as $m\geq 1$.

In the same way we have
$$
\begin{array}{l}
\displaystyle\EE\left[\left(\Xa-P\right)^{2m+1}-\Xa^{2m+1}\right]+~P^{2m+1}\\
\\
\displaystyle=-\sum_{1\leq k\leq m}~\left[W(2(m-k)+2,2k-1)-W(2(m-k)+1,2k)\right]\\
\\
\displaystyle\leq -\frac{1}{2}~\sum_{1\leq k\leq m}~\left(1-\frac{1}{2k}\right)~W(2(m-k)+2,2k-1)\\
\\
\displaystyle\leq -\frac{1}{2}~\left(1-\frac{1}{2m}\right)~\sum_{1\leq k\leq m}~W(2(m-k)+2,2k-1)
\end{array}
$$
This ends the proof of (\ref{estimate-XX-n}).
\section*{Appendix}

\subsection*{Proof of proposition~\ref{letac-theo}}\label{proof-letac-theo}
The lemma is a matrix version of (\ref{Letac-formula}).  As shown in~\cite{letac1} for any sequence of symmetric matrices $Q_i\in \SS_r$
and for any $n\geq 1$
we have
$$
\EE\left(\prod_{1\leq i\leq n}~\langle \Xa,Q_i\rangle_F\right)=\sum_{\sigma\in \Ga_n}~2^{n-\vert \sigma\vert}~\prod_{c\in C(\sigma)}~\tr\left(\prod_{i\in c}(PQ_i)\right)
$$
On the other hand, we have
$$
\left[\Xa\bullet Q\right]^n=\left[\prod_{1\leq i<n}~\langle \Xa,Q_i\rangle_F\right]~\Xa Q_n
$$
 We further assume that
 $$
 Q_n=\frac{1}{2}~\left(e_ke_l^{\prime}+e_le_k^{\prime}\right)
 $$
 for some $1\leq k,l\leq r$, where $e_k$ stands for the (column) vector on the unit sphere 
 defined for any $1\leq i\leq r$ by $e_k(i)=1_{k=i}$. In this case, for any $\sigma\in\Ga_n$ we have
 $$
 \prod_{c\in \sigma}~\tr\left(\prod_{i\in c}(PQ_i)\right)=\left\{
 \prod_{c\in C(n,\sigma)}~\tr\left(\prod_{i\in c}(PQ_i)\right)\right\}~\tr\left(\prod_{i\in c(n,\sigma)}(PQ_i)\right)
 $$
Let $q$ be the length of the cycle $c(n,\sigma)$. Observe that
 $$
 \begin{array}{l}
 \exists 1\leq p\leq q\quad 
c_p(n,\sigma)=n\\
 \\
 \begin{array}[t]{rcl}
 \Longrightarrow  &&
 2\tr\left(\prod_{i\in c(n,\sigma)}(PQ_i)\right)\\
 \\
 &&=2\tr\left[(PQ_{c_1(n,\sigma)})
 \ldots(PQ_{c_{p-1}(n,\sigma)})(PQ_{n})(PQ_{c_{p+1}(n,\sigma)})\ldots (PQ_{c_{q}(n,\sigma)})\right]\\
 &&\\
&&=\tr\left[
Q_{n}(PQ_{c_{p+1}(n,\sigma)})\ldots(PQ_{c_{q}(n,\sigma)})(PQ_{c_1(n,\sigma)})\ldots(PQ_{c_{p-1}(n,\sigma)})P\right]\\
&&\\
 &&=2\left[\left(P\bullet Q\right)^{c^n(\sigma_n)}P\right]_{ sym}(k,l)
 \end{array}
  \end{array}
 $$
 This ends the proof of the first assertion.  The last assertion comes from the fact that
$$
Q_i=P^{v_i}\Longrightarrow
(P\bullet Q)^{c}=P^{\vert c\vert+v(c)}
\quad\mbox{\rm and}\quad \left[\left(P\bullet Q\right)^{c^{\flat}(n,\sigma_n)}P\right]_{ sym}Q_n
=P^{\vert c(n,\sigma_n)\vert +v(c(n,\sigma_n))}
$$
 This ends the proof of the theorem.\cqfd
\subsection*{Proof of the recursion (\ref{recursion-XX-moments-intro})}\label{recursion-XX-moments-intro-proof}

Every permutation $\sigma$ of $\Ga_{m+1}$ can be decomposed
 in an unique way as
$
\sigma=\nu\circ\tau_{i}
$,
for some $\nu\in \Ga_{m}$ (extended to $\Ga_{m+1}$ by setting $\nu(m+1)=m+1$) and a transposition $\tau_{i}$ of the indices
$i$ and $(m+1)$ with $1\leq i\leq m$. When $i=(m+1)$ we clearly have the cycle decomposition
$$
\nu=(c_1(\nu))\ldots (c_p(\nu))\Longrightarrow
\nu\circ\tau_{m+1}=(c_1(\nu))\ldots (c_p(\nu))~(m+1)
$$
This implies that the only cycle containing $(m+1)$ is the $1$-cycle
$$
c(m+1,\nu\circ\tau_{m+1})=(m+1)
$$
and the cycles that doesn't contained $(m+1)$ are given by all the cycles of $\nu$; that is we have
$$
C(m+1,\nu\circ\tau_{m+1})=C(\nu)
$$
This yields the decomposition
$$
\begin{array}{l}
\displaystyle\sum_{\sigma\in\Ga_{m+1}}~\left(\frac{1}{2}\right)^{\vert \sigma\vert}
 \tr_{C(m+1,\sigma)}(P^{\overline{v}})~
 P^{\overline{v}(c(m+1,\sigma))}\\
 \\
 =\displaystyle \frac{1}{2}~P^{1+v_{m+1}}~\sum_{\nu\in\Ga_{m}}~\left(\frac{1}{2}\right)^{\vert \nu\vert}
\tr_{C(\nu)}(P^{\overline{v}})~
 \displaystyle+\sum_{1\leq i\leq m}\sum_{\nu\in\Ga_{m}}~\left(\frac{1}{2}\right)^{\vert \nu\circ\tau_{i}\vert}
\tr_{C(m+1,\nu\circ\tau_{i})}(P^{\overline{v}})~
 P^{\overline{v}(c(m+1,\nu\circ\tau_{i}))}
\end{array}
$$
On the other hand, for any $1\leq i\leq m$ the cycle $c(m+1,\nu\circ\tau_{i})$ can be expressed as
$$
\left(c_1\rightarrow \ldots \rightarrow c_k=i\rightarrow m+1\rightarrow c_{k+1}\rightarrow\ldots \rightarrow c_q\rightarrow c_1\right)
$$
where
$$
\left(c_1\rightarrow \ldots \rightarrow c_k=i\rightarrow c_{k+1}\rightarrow\ldots \rightarrow c_q\rightarrow c_1\right)=c(i,\nu)
$$
is the cycle of $\nu$ containing $i$. This shows that
$$
\vert c(m+1,\nu\circ\tau_{i})\vert=\vert c(i,\nu)\vert+1
\Longrightarrow
\overline{v}( c(m+1,\nu\circ\tau_{i}))=
\vert c(i,\nu)\vert+1+v(c(i,\nu))+v_{m+1}
$$
and
$$
\vert \nu\circ\tau_{i}\vert=\vert C(i,\nu) \vert+1=\vert \nu\vert
$$
This yields the formula
$$
\begin{array}{l}
\displaystyle\sum_{\sigma\in\Ga_{m+1}}~\left(\frac{1}{2}\right)^{\vert \sigma\vert}
 \tr_{C(m+1,\sigma)}(P^{\overline{v}})~
 P^{\overline{v}(c(m+1,\sigma))}\\
 \\
 =\displaystyle ~P^{1+v_{m+1}}~\left[\frac{1}{2}~\sum_{\nu\in\Ga_{m}}~\left(\frac{1}{2}\right)^{\vert \nu\vert}
\tr_{C(\nu)}(P^{\overline{v}})~
 \displaystyle+\sum_{1\leq i\leq m}\sum_{\nu\in\Ga_{m}}~\left(\frac{1}{2}\right)^{\vert \nu \vert}
\tr_{C(i,\nu)}(P^{\overline{v}})~
 P^{\overline{v}(c(i,\nu))}\right]
\end{array}
$$
The end of the proof of (\ref{recursion-XX-moments-intro}).\cqfd

\subsection*{Proof of lemma~\ref{lem-tech}}\label{proof-lem-tech}

We have
\begin{eqnarray*}
f(n)&\leq & \frac{1}{n}~g(0)~\rho(g)~f(n-1)\\
&&\hskip1cm+\left(1-\frac{1}{n}\right)~\frac{1}{n-1}~\sum_{0\leq k<(n-1)}~\frac{g(1+((n-1)-k))}{g((n-1)-k)}~g((n-1)-k)~f(k)\\
&\leq &\rho(g)~\left[\frac{1}{n}~g(0)+\left(1-\frac{1}{n}\right)~\right]~f(n-1)\leq \rho(g)^{n-m}~\prod_{m<l\leq n}\left[1-\frac{1}{l}\left(1-g(0)\right)\right]~ f(m)
\end{eqnarray*}
This implies that
$$
f(n)/f(m)\leq  \rho(g)^{n-m}~\prod_{m< l\leq n}\left[1-\frac{1-g(0)}{l}
\right]=(g(0)\rho(g))^{n-m}~\frac{m!}{n!}~\prod_{m\leq l<n}(1+l/g(0))
$$
When $g(0)\leq 1/2$
\begin{eqnarray*}
\prod_{m\leq l<n}(1+l/g(0))&\leq &\prod_{m\leq l<n}(1+2l)= (2m+1)(2m+3)\ldots (2n-1)\\
&=&2^{-(n-m)}\frac{(2m+1)(2m+2)\ldots (2n-1)(2n)}{(m+1)\ldots n}=2^{-(n-m)}~\frac{(2n)!}{(2m)!}~\frac{m!}{n!}
\end{eqnarray*}
This implies that
$$
f(n)\leq \rho(g)^{n-m}~\frac{2^{-2n}(2n)!/n!^2}{2^{-2m}(2m)!/m!^2}~
~f(m)
$$
By Stirling approximation 
$$
\exp{\left(\frac{1}{12n+1}\right)}\leq \frac{n!}{\sqrt{2\pi n}~n^n~e^{-n}}\leq \exp{\left(\frac{1}{12n}\right)}
$$
we have
\begin{eqnarray*}
2^{-2n}~\frac{(2n)!}{n!^2}&\leq&~\frac{1}{\sqrt{\pi n}}~\exp{\left(\frac{1}{6n}\left[\frac{1}{4}-\frac{6n}{6 n+1/2}\right]\right)}\\
&=&~\frac{1}{\sqrt{\pi n}}~\exp{\left(-\frac{1}{6n}\left[\frac{3}{4}-\frac{1}{12 n+1}\right]\right)}\leq\frac{1}{\sqrt{\pi n}}~\exp{\left(-\frac{1}{12n}\left[1-\frac{1}{9 n}\right]\right)}
\end{eqnarray*}
This yields
$$
f(n)\leq (\rho(g))^{n}~\frac{1}{\sqrt{\pi n}}~\exp{\left(-\frac{1}{12n}\left[1-\frac{1}{9 n}\right]\right)}
~f(0)
$$
This ends the proof of the lemma.\cqfd

 \subsection*{Some combinatorial formulae}
We recall that number of solutions to the equation
$$
j_1+j_2+\ldots+j_{n-k}=k-l
$$
with $j_1,\ldots,j_{n-k}\in\NN$ is given by the binomial coefficient
$$
\left(
\begin{array}{c}
(k-l)+(n-k)-1\\
(k-l)
\end{array}
\right)=\left(
\begin{array}{c}
n-(l+1)\\
k-l
\end{array}
\right)
$$
In addition we have
\begin{lem}\label{bino-appendix}
$$
\left(
\begin{array}{c}
n+1\\
k
\end{array}
\right)=\sum_{0\leq l\leq k}~(l+1) \left(
\begin{array}{c}
n-(l+1)\\
k-l
\end{array}
\right)
$$ 
\end{lem}
\proof
We use Pascal triangle formula
$$
\left(
\begin{array}{c}
n+1\\
k
\end{array}
\right)=\left(
\begin{array}{c}
n\\
k
\end{array}
\right)+\left(
\begin{array}{c}
n\\
k-1
\end{array}
\right)
$$
This readily implies that
\begin{eqnarray*}
k~\left(
\begin{array}{c}
n-k+1\\
1
\end{array}
\right)+\left(
\begin{array}{c}
n-k+1\\
0
\end{array}
\right)&=&k~\left[\left(
\begin{array}{c}
n-k\\
0
\end{array}
\right)+\left(
\begin{array}{c}
n-k\\
1
\end{array}
\right)\right]+\left(
\begin{array}{c}
n-k+1\\
0
\end{array}
\right)\\
&=&k~\left(
\begin{array}{c}
n-k\\
0
\end{array}
\right)+(k+1)\\
&=&k\left(
\begin{array}{c}
n-k\\
0
\end{array}
\right)+(k+1)~\left(
\begin{array}{c}
n-(k+1)\\
0
\end{array}
\right)
\end{eqnarray*}
 This shows that
 $$
 \begin{array}{l}
\sum_{0\leq l\leq k}~(l+1) \left(
\begin{array}{c}
n-(l+1)\\
k-l
\end{array}
\right)\\
\\
=\sum_{0\leq l\leq k-2}~(l+1) \left(
\begin{array}{c}
n-(l+1)\\
k-l
\end{array}
\right)+k~\left(
\begin{array}{c}
n-k+1\\
1
\end{array}
\right)+\left(
\begin{array}{c}
n-k+1\\
0
\end{array}
\right)
\end{array}
 $$
Now we prove the following formula
 $$
 \begin{array}{l}
 \left(
\begin{array}{c}
n+k\\
k
\end{array}
\right)\\
\\
=\sum_{0\leq l\leq m}~(l+1) \left(
\begin{array}{c}
n-(l+1)\\
k-l
\end{array}
\right)+(m+2)~\left(
\begin{array}{c}
n-(m+1)\\
k-(m+1)
\end{array}
\right)+\left(
\begin{array}{c}
n-(m+1)\\
k-(m+2)
\end{array}
\right)
\end{array}
 $$
 by induction w.r.t. the parameter $0\leq m\leq k-2$. The result is immediate for $m=0$ using Pascal triangle formula. Indeed we have 
  \begin{eqnarray*}
 \left(
\begin{array}{c}
n+1\\
k
\end{array}
\right)
&=& \left(
\begin{array}{c}
n\\
k
\end{array}
\right) + \left(
\begin{array}{c}
n\\
k-1
\end{array}
\right)\\
&=&\left(
\begin{array}{c}
n-1\\
k
\end{array}
\right) +\left(
\begin{array}{c}
n-1\\
k-1
\end{array}
\right) 
+ \left(
\begin{array}{c}
n-1\\
k-1
\end{array}
\right)+\left(
\begin{array}{c}
n-1\\
k-2
\end{array}
\right)\\
&=&\left(
\begin{array}{c}
n-1\\
k
\end{array}
\right)+2~\left(
\begin{array}{c}
n-1\\
k-1
\end{array}
\right)+\left(
\begin{array}{c}
n-1\\
k-2
\end{array}
\right)
\end{eqnarray*}
 We further assume that the formula is true at rank $m$. In this situation, using Pascal formula we have
   $$
  \begin{array}{l}
 (m+2)~\left(
\begin{array}{c}
n-(m+1)\\
k-(m+1)
\end{array}
\right)+\left(
\begin{array}{c}
n-(m+1)\\
k-(m+2)
\end{array}
\right)\\
\\
= (m+2)~\left[
\left(
\begin{array}{c}
n-(m+2)\\
k-(m+1)
\end{array}
\right)
+
\left(
\begin{array}{c}
n-(m+2)\\
k-(m+2)
\end{array}
\right)
\right]\\
\\
\hskip3cm+\left[\left(
\begin{array}{c}
n-(m+2)\\
k-(m+2)
\end{array}
\right)+[\left(
\begin{array}{c}
n-(m+2)\\
k-(m+3)
\end{array}
\right)\right]\\
\\
=(1+(m+1))~
\left(
\begin{array}{c}
n-((m+1)+1)\\
k-(m+1)
\end{array}
\right)+(m+3)~
\left(
\begin{array}{c}
n-(m+2)\\
k-(m+2)
\end{array}
\right)+\left(
\begin{array}{c}
n-(m+2)\\
k-(m+3)
\end{array}
\right)
\end{array}
 $$
 The end of the proof of the induction is now clear, thus it is skipped.
 This completes the proof of the lemma.
 \cqfd
 
 Now we come to the proof of (\ref{C-n-m-n}).
 
 {\bf Proof of (\ref{C-n-m-n}):}\label{proof-C-n-m-n}
 The above lemma implies that
\begin{eqnarray*}
\displaystyle \sum_{v_1+\ldots+v_m=n}(v_{m}+1)
&=&\sum_{v_m=0}^n  (v_{m}+1)~\sum_{v_1+v_{2}+\ldots+v_{m-1}=n-v_m} 1\\
&=&\sum_{v_m=0}^n  (v_{m}+1)~
   \left(
  \begin{array}{c}
  (n-v_m)+(m-1)-1\\
  (n-v_m)
  \end{array}
  \right)\\
& =&\sum_{v_m=0}^n  (v_{m}+1)~
   \left(
  \begin{array}{c}
  (n+(m-1))-(v_m+1)\\
  n-v_m
  \end{array}
  \right)= \left(
  \begin{array}{c}
  n+m\\
  n
  \end{array}
  \right)
\end{eqnarray*}

 \subsection*{Proof of lemma~\ref{F-lemma}}\label{proof-F-lemma}
 
Observe that
\begin{eqnarray*}
\displaystyle V_{m,n+1}
\displaystyle &\subset&\cup_{1\leq i\leq m}\{(v_1,\ldots,v_{m})\in V_{m,n+1}~\mbox{\rm 
s.t.}\quad v_i\geq 1\}=\cup_{1\leq i\leq m}\{\theta_i(v,0)~:~
v\in V_{m,n}\}
 \end{eqnarray*}
with the mappings $\theta_i$ introduced in definition~\ref{def-theta-i}. We recall that 
$$
\forall  (v,w)\in (V_{m,n}\times \NN)
\quad \theta_i(v,w)=(v_1,\ldots,v_{i-1},1+v_i+w,v_{i+1},\ldots,v_{m})
$$

This yields for any non negative function $F$ the estimate
\begin{eqnarray*}
\displaystyle 
\sum_{v\in V_{m,n+1}}~F(v)
& \leq &\sum_{1\leq i\leq m}~\sum_{v\in V_{m,n}}~
F(\theta_i(v,0))\\
&=&\displaystyle \sum_{1\leq i\leq m}~\sum_{v_1+\ldots+v_{m}+0=n}~F(v_1,\ldots,v_{i-1},v_i+1+0,v_{i+1},\ldots,v_{m})\\
&\leq& \sum_{1\leq j\leq m}~\sum_{v\in V_{m+1,n}}~F\left(\theta_j(v)\right)\leq \frac{1}{2~\left(1+\epsilon_m\right)}~\sum_{v\in V_{m+1,n}}~F(v)
 \end{eqnarray*}
from which we find the estimate (\ref{first-estimate}).

We want to check that
$$
\begin{array}{l}
\displaystyle\sum_{ v\in V_{m,n+1}}~(\pi^+(v)+1)~F(v)\leq
 \sum_{ v\in V_{m+1,n}}~(\pi^+(v)+1)~F(v)
\end{array}
$$
By (\ref{Hyp-F}) we have
$$
\begin{array}{l}
 \displaystyle \sum_{ v\in V_{m+1,n}}~(\pi^+(v)+1)~F(v)\\
 \\
 \displaystyle\geq 2~(1+\epsilon_m)~\sum_{1\leq j\leq m}~ 
 \sum_{v\in V_{m+1,n}}~(\pi^+(v)+1)~
 F\left(\theta_j(v)\right)\\
 \\
 = \displaystyle 2~(1+\epsilon_m)
 \sum_{1\leq j\leq m}~\sum_{0\leq \vert v\vert_{-j}\leq n}\\
 \\
\hskip1cm \displaystyle \left[\sum_{v_j+v_{m+1}= {n}-\vert v\vert_{-j}}~(v_{m+1}+1)~\right]~F\left(v_1,\ldots,v_{j-1},1+\left(n-\vert v\vert_{-j}\right),v_{j+1},\ldots,v_m\right)
\end{array}
$$
with
$$
\vert v\vert_{-j}:=\sum_{i\in [m]-\{j\}}v_i
$$
By (\ref{C-n-m-n}) this implies that
$$
\begin{array}{l}
 \displaystyle \sum_{ v\in V_{m+1,n}}~(\pi^+(v)+1)~F(v)\\
 \\
 \displaystyle\geq  \displaystyle 2~(1+\epsilon_m)
 \sum_{1\leq j\leq m}~\sum_{0\leq \vert v\vert_{-j}\leq n} \left(
\begin{array}{c}
 {n+2}-\vert v\vert_{-j}\\
 {n}-\vert v\vert_{-j}
\end{array}
\right)~F\left(v_1,\ldots,v_{j-1},1+\left(n-\vert v\vert_{-j}\right),v_{j+1},\ldots,v_m\right)
\end{array}
$$
Changing the index $1+ {n}-\vert v\vert_{-j}$ by $v_j\geq 1$ in each summand we find that
\begin{equation}\label{summands-refs}
\begin{array}{l}
 \displaystyle\sum_{ v\in V_{m+1,n}}~(\pi^+(v)+1)~F(v)\\
 \\
 \displaystyle\geq  \displaystyle (1+\epsilon_m)
 \sum_{1\leq j\leq m}~\sum_{v_1+\ldots+v_m=n+1,~v_j\geq 1}~ v_j\left(1+v_j\right)
~F\left(v_1,\ldots,v_m\right)\\
 \\
 \displaystyle=  \displaystyle (1+\epsilon_m)
~\sum_{v_1+\ldots+v_m=n+1}~ \left[ \sum_{1\leq j\leq m}v_j\left(1+v_j\right)\right]
~F\left(v_1,\ldots,v_m\right)
\end{array}
\end{equation}
On the other hand we have
$$
\begin{array}{l}
\displaystyle\sum_{v_1+\ldots+v_{m}= {n+1}}~(v_{m}+1)~F(v_1,\ldots,v_{m})\\
\\
=\displaystyle\sum_{v_1+\ldots+v_{m-1}= {n+1}}~~F(v_1,\ldots,v_{m-1},0)+\sum_{v_1+\ldots+v_{m}= {n+1},~v_m\geq 1}~(v_{m}+1)~F(v_1,\ldots,v_{m})
\end{array}$$
The summand in (\ref{summands-refs}) associated with $j=m$ is given by
$$
\sum_{v_1+\ldots+v_m=n+1,~v_m\geq 1}~ v_m\left(1+v_m\right)
~F\left(v_1,\ldots,v_m\right)\geq \sum_{v_1+\ldots+v_{m}= {n+1},~v_m\geq 1}~(v_{m}+1)~F(v_1,\ldots,v_{m})
$$
We also have
$$
\begin{array}{l}
\displaystyle\sum_{v_1+\ldots+v_{m-1}= {n+1}}~~F(v_1,\ldots,v_{m-1},0)\\
\\
\leq 
\displaystyle\sum_{1\leq j\leq m-1 }\sum_{v_1+\ldots+v_{m-1}+0= {n+1},~v_j\geq 1}~~F(v_1,\ldots,v_{m-1},0)\\
\\
\displaystyle\leq  \sum_{1\leq j\leq m-1}~\sum_{v_1+\ldots+v_{m-1}+v_m=n+1,~v_j\geq 1}~ v_j\left(1+v_j\right)
~F\left(v_1,\ldots,v_m\right)
\end{array}$$
This ends the proof of the lemma.
\cqfd

\section*{Proof of the estimates (\ref{case-PI-intro-estimate-r1}) and (\ref{case-PI-intro-estimate-r2})}\label{proof-PI-intro-estimate-r12}

The binomial formula (\ref{case-PI-intro}) is a direct consequence of (\ref{case-PI-2}).

\subsection*{The case $r=1$}
We have
\begin{eqnarray*}
\EE(X^{2n})&=&\frac{(2n)!}{n!}~2^{-n}=(2n-1)~\EE(X^{2(n-1)})
\end{eqnarray*}
This shows that
$$
\EE((X^2-1)^n)=z_n(0)+\sum_{1\leq k\leq n}~(-1)^k~z_n(k)
$$
with 
\begin{eqnarray*}
z_n(k)&:=&{n \choose k}~\EE(X^{2(n-k)})=(2(n-k)-1)~\frac{k+1}{n-k}~{n \choose k+1}~\EE(X^{2(n-(k+1))})\\
&=&(k+1)~\left(2-\frac{1}{n-k}\right)~z_n(k+1) ~> (k+1)~z_n(k+1) 
\end{eqnarray*}
for any $k<n$.
Equivalently, we have
$$
z_n(k+1) =  \frac{1}{2(k+1)}\frac{1}{\left(1-\frac{1}{2(n-k)}\right)}~z_n(k)
$$
Assume that $n=2m$ is an even integer. In this case
$$
\begin{array}{l}
\displaystyle\sum_{1\leq k\leq n}~(-1)^k~z_n(k)\\
\\
\displaystyle=\sum_{1\leq k\leq m}~\left(z_{2m}(2k)-z_{2m}(2k-1)\right)=-
\sum_{1\leq k\leq m}~\left(z_{2m}(2k-1)-z_{2m}((2k-1)+1)\right)\\
\\
\displaystyle= -\frac{1}{2}~\sum_{1\leq k\leq m}~\left[2-\frac{1}{(k+1)~\left(1-\frac{1}{2(2m-k)}\right)}
\right]~z_{2m}(2k-1)\\
\\
\displaystyle=-\frac{1}{2}~\sum_{1\leq k\leq m}~\left(1+\frac{1}{2(2m-k)-1}\right)~ \left[
2\left(1-\frac{1}{2(2m-k)}\right)-\frac{1}{k+1}
\right]~z_{2m}(2k-1)\\
\\
\displaystyle=-\frac{1}{2}~\left(1+\frac{1}{2(2m-1)-1}\right)~\sum_{1\leq k\leq m} \left[
\left(1-\frac{1}{2m-k}\right)+\left(1-\frac{1}{k+1}\right)
\right]~z_{2m}(2k-1)
\end{array}
$$
This implies that
\begin{eqnarray*}
\sum_{1\leq k\leq n}~(-1)^k~z_n(k)&\leq& -\frac{2m-1}{4m-3} \left[
\left(1-\frac{1}{m}\right)+\frac{1}{2}
\right]~\sum_{1\leq k\leq m}~z_{2m}(2k-1)\\
&=& -\frac{2m-1}{4m-3}~\frac{3m-2}{2m}~\sum_{1\leq k\leq m}~z_{2m}(2k-1)\\
&=&  -  \frac{3}{4} ~\frac{m-1/2}{m-3/4}~\left(1-\frac{2}{3m}\right)~\sum_{0\leq k<m}~{2m \choose 2k+1}~\EE(X^{2(2k+1)})
\end{eqnarray*}
the last assertion comes from the fact that
\begin{eqnarray*}
\sum_{1\leq k\leq m}~z_{2m}(2k-1)&=&\sum_{1\leq k\leq m}~
{2m \choose 2(m-k)+1}~\EE(X^{2(2(m-k)+1)})\\
&=&\sum_{0\leq k<m}~{2m \choose 2k+1}~\EE(X^{2(2k+1)})
\end{eqnarray*}
This ends the proof of the  estimate (\ref{case-PI-intro-estimate-r1}).

When $n=2m+1$ is an odd integer, arguing as above  we have
$$
\begin{array}{l}
\displaystyle
\sum_{1\leq k\leq n}~(-1)^k~z_n(k)\\
\\
\displaystyle=-z_{2m+1}(2m+1)+
\sum_{1\leq k\leq m}~\left(z_{2m+1}(2k)-z_{2m+1}(2k-1)\right)\\
\\\displaystyle
\leq -1-\frac{1}{2}~\sum_{1\leq k\leq m}~\left[2-\frac{1}{(k+1)~\left(1-\frac{1}{2(2m+1-k)}\right)}
\right]~z_{2m+1}(2k-1)\\
\\
\displaystyle
\leq -1-\frac{1}{2}~\sum_{1\leq k\leq m}~\left[1+\frac{1}{2(2m+1-k)-1}\right]~\left[2\left(1-\frac{1}{2(2m+1-k)}\right)-\frac{1}{k+1}
\right]~z_{2m+1}(2k-1)\\
\end{array}
$$
This implies that
\begin{eqnarray*}
\displaystyle
\sum_{1\leq k\leq n}~(-1)^k~z_n(k)
&\leq& -1-\frac{m}{4m-1}~\frac{3m+1}{m+1}~\sum_{1\leq k\leq m}
~z_{2m+1}(2k-1)\\
&
\leq& -1-\frac{m}{4m-1}~\frac{3m+1}{m+1}~\sum_{1\leq k\leq m}
~{2m+1 \choose 2((m-k)+1)}~\EE(X^{4((m-k)+1)})\\
&=& -1-\frac{m}{4m-1}~\frac{3m+1}{m+1}~\sum_{1\leq k\leq m}
~{2m+1 \choose 2k}~\EE(X^{4k})
\end{eqnarray*}
as soon as $m\geq 1$.  This yields the estimate
$$
\EE((X^2-1)^{2m+1})\leq \EE(X^{4m+2}) -1-\frac{m}{4m-1}~\frac{3m+1}{m+1}~\sum_{1\leq k\leq m}
~{2m+1 \choose 2k}~\EE(X^{4k})$$

\subsection*{The case $r\geq 2$}

Let $R=\langle X,X\rangle$ and $\Ua=\Xa/R$. Recalling that $R$ and $\Ua$ are independent and using
 (\ref{case-PI-2}) we readily check that
\begin{eqnarray*}
\EE\left[\left(\Xa-I\right)^{n}\right]
 &  =&\sum_{0\leq  k\leq n}
\left(
\begin{array}{c}
n\\
k
\end{array}
\right)~(-1)^{n-k}~\EE\left[ \Xa^{k}\right]
 \displaystyle  =\sum_{0\leq  k\leq n}
\left(
\begin{array}{c}
n\\
k
\end{array}
\right)~(-1)^{n-k}~\EE\left[ R^k~\Ua\right]\\
 & =&
 \sum_{0\leq  k\leq n}
\left(
\begin{array}{c}
n\\
k
\end{array}
\right)~(-1)^{n-k}~\EE\left[ R^k\right]~\frac{1}{r}~I\\
&=& \displaystyle  \frac{1}{r}
 \sum_{0\leq  k\leq n}
\left(
\begin{array}{c}
n\\
k
\end{array}
\right)~(-1)^{n-k}~\left[\prod_{0\leq l<k}\left(r+2l\right)\right]~I
\end{eqnarray*}
 This implies that
\begin{eqnarray*}
\EE\left[\left(\Xa-I\right)^{n}\right]
&=& \displaystyle  z_{n}(0)~I+
 \sum_{1\leq  k\leq n}~(-1)^{k}~z_n(k)~I\\
& =& \displaystyle  \EE\left[\Xa^{n}\right]+
 \sum_{1\leq  k\leq n}~(-1)^{k}~z_n(k)~I
\end{eqnarray*}
with the collection of parameters
\begin{eqnarray*}
z_n(k)&:=&\left(
\begin{array}{c}
n\\
k
\end{array}
\right)~\frac{1}{r}~\prod_{0\leq l<n-k}\left(r+2l\right)\\
&=&\frac{k+1}{n-k}~(2(n-k)-(2-r))~z_n(k+1)=2(k+1)~\left(1-\frac{1-r/2}{n-k}\right)~z_n(k+1)
\end{eqnarray*}
This implies that
$$
z_n(k+1)= 
 \frac{1}{2(k+1)\left(1+\frac{r-2}{2(n-k)}\right)}~ z_n(k)
$$
Assume that $n=2m$ is an even integer and 
 $r\geq 2$. In this case we have
$$
\begin{array}{l}
\displaystyle\sum_{1\leq k\leq n}~(-1)^k~z_n(k)\\
\\
\displaystyle=\sum_{1\leq k\leq m}~\left(z_{2m}(2k)-z_{2m}(2k-1)\right)=-
\sum_{1\leq k\leq m}~\left(z_{2m}(2k-1)-z_{2m}((2k-1)+1)\right)\\
\\
\displaystyle= -\frac{1}{2}~\sum_{1\leq k\leq m}~\left[2-\frac{1}{(k+1)~\left(1+\frac{r/2-1}{(2m-k)}\right)}
\right]~z_{2m}(2k-1)\\
\\
\displaystyle=-\frac{1}{2}~\sum_{1\leq k\leq m}~\left(1-\frac{r/2-1}{(2m-k)+(r/2-1)}\right)~ \left[
2\left(1+\frac{r/2-1}{(2m-k)}\right)-\frac{1}{k+1}
\right]~z_{2m}(2k-1)\\
\\
\displaystyle\leq -\frac{1}{2}~\left(1-\frac{r/2-1}{m+(r/2-1)}\right)~\sum_{1\leq k\leq m} \left[
\left(1+\frac{r-2}{2m-k}\right)+\left(1-\frac{1}{k+1}\right)
\right]~z_{2m}(2k-1)\\
\\
\displaystyle\leq -\frac{1}{2}~\left(1-\frac{r/2-1}{m+(r/2-1)}\right)~\left(\frac{3}{2}+\frac{r-2}{2m-1}\right)~\sum_{1\leq k\leq m}~z_{2m}(2k-1)
\end{array}
$$
This implies that
$$
\begin{array}{l}
\displaystyle\sum_{1\leq k\leq n}~(-1)^k~z_n(k)\\
\\
\displaystyle\leq-\frac{m}{m+(r/2-1)}~\left(\frac{3}{4}+\frac{r/2-1}{2m-1}\right)
\sum_{1\leq k\leq m}~z_{2m}(2k-1)
\end{array}
$$
We conclude that
$$
\begin{array}{l}
   \displaystyle\EE\left[\left(\Xa-I\right)^{2m}\right]
\\
\\
\leq \displaystyle  \EE\left[\Xa^{2m}\right]-
\frac{m}{2m+(r-2)}~\left(\frac{3}{2}+\frac{r-2}{2m-1}\right)
\sum_{1\leq k\leq m}~z_{2m}(2k-1)~I\\
\\
= \displaystyle  \EE\left[\Xa^{2m}\right]-
\frac{m}{2m+(r-2)}~\left(\frac{3}{2}+\frac{r-2}{2m-1}\right)
\sum_{1\leq k\leq m}~
\left(
\begin{array}{c}
2m\\
2k-1
\end{array}
\right)~\EE\left[\Xa^{2k-1}\right]~I
\end{array}
$$
This end the proof of (\ref{case-PI-intro-estimate-r2}).

When $n=2m+1$ is an odd integer and $r\geq 2$, arguing as above  we have
$$
\begin{array}{l}
\displaystyle
\sum_{1\leq k\leq n}~(-1)^k~z_n(k)\\
\\
\displaystyle=-z_{2m+1}(2m+1)-
\sum_{1\leq k\leq m}~\left(z_{2m+1}(2k-1)-z_{2m+1}(2k)\right)\\
\\\displaystyle
\leq -\frac{1}{r}-\frac{1}{2}~\sum_{1\leq k\leq m}~\left[2-\frac{1}{(k+1)~\left(1+\frac{r-2}{2(2m+1-k)}\right)}
\right]~z_{2m+1}(2k-1)\\
\\
\displaystyle
= -\frac{1}{r}-\frac{1}{2}~\sum_{1\leq k\leq m}~\left[1-\frac{r-2}{2(2m+1-k)+(r-2)}\right]~\left[2\left(1+\frac{r-2}{2(2m+1-k)}\right)-\frac{1}{k+1}
\right]~z_{2m+1}(2k-1)\\
\\
\displaystyle
\leq -\frac{1}{r}-\frac{1}{2}~\left[1-\frac{r-2}{2m+r}\right]~\left(\frac{3}{2}+\frac{r-2}{2m}\right)
~\sum_{1\leq k\leq m}~z_{2m+1}(2k-1)\
\end{array}
$$
This yields the estimate
$$
\begin{array}{l}
\displaystyle
\sum_{1\leq k\leq n}~(-1)^k~z_n(k)\\
\\
\displaystyle
\leq -\frac{1}{r}-\frac{1}{4}~\frac{m+1}{m+r/2}~\left(3+\frac{r-2}{m}\right)
~\sum_{1\leq k\leq m}~z_{2m+1}(2k-1)\\
\\
\displaystyle
\leq -\frac{1}{r}-\frac{1}{4}~\frac{m+1}{m+r/2}~\left(3+\frac{r-2}{m}\right)
~\sum_{1\leq k\leq m}~\left(
\begin{array}{c}
2m+1\\
2k-1
\end{array}
\right)~\EE\left[\Xa^{2m+1-(2k-1)}\right]\\
\end{array}
$$
from which we conclude that
$$
\begin{array}{l}
\displaystyle\EE\left[\left(\Xa-I\right)^{2m+1}\right]\\
\\
\leq \displaystyle    \EE\left[\Xa^{2m+1}\right] -\frac{1}{r}-\frac{1}{4}~\frac{m+1}{m+r/2}~\left(3+\frac{r-2}{m}\right)
\sum_{1\leq k\leq  m}~\left(
\begin{array}{c}
2m+1\\
2k
\end{array}
\right)~\EE\left[\Xa^{2k}\right]
\end{array}
$$
\subsection{Proof of the Laplace formula (\ref{t-X-Laplace}) and (\ref{alternative-recursion})}\label{proof-estimate-Laplace-xx}

By definition of the  Wishart distribution
with a single degree of freedom and covariance matrix $P$, we have
$$
\begin{array}{l}
\EE\left(\Xa~e^{t\Xa}\right)\\
\\
=\EE\left(\Xa~e^{t\langle X,X\rangle}\right)\\
\\
\displaystyle=\int~\frac{1}{2\Gamma_r(1/2)~\mbox{\rm det}(P)^{1/2}}~
\mbox{\rm det}(Q)^{(1-r-1)/2}~Q~e^{t\,\tr(Q)}~\exp{\left[-\frac{1}{2}~\tr\left(P^{-1}Q\right)\right]}~\gamma(dQ)\\
\\
\displaystyle=\int~\frac{1}{2\Gamma_r(1/2)~\mbox{\rm det}(P)^{1/2}}~
\mbox{\rm det}(Q)^{(1-r-1)/2}~Q~\exp{\left[-\frac{1}{2}~\tr\left([I-2t~P]P^{-1}Q\right)\right]}~\gamma(dQ)\\
\\
\displaystyle=\int~\frac{1}{2\Gamma_r(1/2)~\mbox{\rm det}(P)^{1/2}}~
\mbox{\rm det}(Q)^{(1-r-1)/2}~Q~\exp{\left[-\frac{1}{2}~\tr\left(\left\{[I-2t~P]^{-1}P\right\}^{-1}Q\right)\right]}~\gamma(dQ)\\
\\
\displaystyle=~\left\{\frac{\mbox{\rm det}([I-2t~P]^{-1}P)}{\mbox{\rm det}(P)}\right\}^{1/2}~[I-2t~P]^{-1}P=
\left\{\frac{1}{\mbox{\rm det}([I-2t~P])}\right\}^{1/2}~[I-2t~P]^{-1}P
\end{array}
$$
where $\gamma(dQ)$ stands for the Lebesgue measure on the cone of symmetric
positive definite matrices, and $\Gamma_{r}$ is the multivariate gamma function
$$
\Gamma_{r}(z)=\pi^{r(r-1)/4}~\prod_{1\leq k\leq r}\Gamma\left(z-\frac{k-1}{2}\right)
$$
This implies that
$$
\begin{array}{l}
\displaystyle\EE\left(\exp{\left(t\,\langle X,X\rangle\right)}\right)\\
\\
=\displaystyle\exp{\left[-\frac{1}{2}\log\left(\mbox{\rm det}([I-2t~P])\right)\right]}=\exp{\left[\frac{1}{2}\sum_{n\geq 1}~\frac{(2t)^n}{n}~\tr(P^n)\right]}~\\
\\
\displaystyle\leq \exp{\left[\frac{1}{2}\sum_{n\geq 1}~\frac{(2t)^n}{n}~\tr(P)^n\right]}=
 \exp{\left[-\frac{1}{2}\log{\left(1-2\tr(P)\right)}\right]}=\frac{1}{\sqrt{1-2t\, \tr(P)}}~
 \end{array}
$$
as well as
\begin{eqnarray*}
\displaystyle\partial_t\EE\left(\exp{\left(t\Xa\right)}\right)
&=&\exp{\left[-\frac{1}{2}\log\left(\mbox{\rm det}([I-2t~P])\right)\right]}~~[I-2t~P]^{-1}P
\end{eqnarray*}
This ends the proofs of (\ref{t-X-Laplace}) and (\ref{alternative-recursion}).\cqfd

\section*{Proof of the polarization formula (\ref{polarization})}\label{proof-polarization}

To check this claim, observe that
\begin{eqnarray*}
\EE\left(\prod_{1\leq i\leq m}\langle X,Q_i X\rangle \right)=\frac{1}{m!}~\sum_{\sigma\in \Ga_m}~
\EE\left(\prod_{1\leq i\leq m}\langle X,Q_{\sigma(i)} X\rangle \right)
\end{eqnarray*}
On the other hand, we have
$$
\begin{array}{l}
\displaystyle\EE\left(\left[\prod_{1\leq i\leq m}W_i\right]~\left(\sum_{1\leq i\leq m}W_i~\langle X,Q_iX\rangle \right]^n\right)\\
\\
=\displaystyle\sum_{a\in [m]^[m]}\EE\left\{\EE\left(
\left[\prod_{1\leq i\leq m}W_i\right]~\left[\prod_{1\leq k\leq m}W_{a(k)}\right]~|~X\right)
\prod_{1\leq k\leq m}\langle X,Q_{a(k)}X\rangle~|~(W_i)_i\right\}\\
\\
=\displaystyle\sum_{\sigma\in\Ga_m}\EE\left\{\EE\left(
\left[\prod_{1\leq i\leq m}W_i\right]~\left[\prod_{1\leq k\leq m}W_{\sigma(k)}\right]~|~X\right)
\prod_{1\leq k\leq m}\langle X,Q_{\sigma(k)}X\rangle\right\}\\
\\
\displaystyle=\sum_{\sigma\in \Ga_m}~
\EE\left(\prod_{1\leq i\leq m}\langle X,Q_{\sigma(i)} X\rangle \right)
\end{array}
$$
The last couple of assertions comes from the fact that
$$
\forall \sigma\in\Ga_n\quad
\left(1,\ldots,n\right)\not=\left(j_{\sigma(1)},\ldots,j_{\sigma(n)}\right)
\Longrightarrow
\EE\left(\left[\prod_{1\leq i\leq n}W_i\right]~\left[\prod_{1\leq k\leq n}W_{j_k}\right]~|~X\right)=0
$$
and
$$
\left(1,\ldots,n\right)=\left(j_{\sigma(1)},\ldots,j_{\sigma(n)}\right)
\Longrightarrow
\EE\left(\left[\prod_{1\leq i\leq n}W_i\right]~\left[\prod_{1\leq k\leq n}W_{j_k}\right]~|~X\right)=1
$$

\end{document}